\documentclass[10pt]{amsart}
\usepackage{amssymb}
\usepackage{MnSymbol}

\usepackage{amscd}
\usepackage{amsthm}
\usepackage[dvips]{graphicx}
\usepackage{color}
\usepackage[all]{xy}
\usepackage{verbatim}

\newtheorem{Lemma1}{{Lemma}}[section]
\newtheorem{Theo1}[Lemma1]{{Theorem}}
\newtheorem*{Theo2}{{Theorem}}
\newtheorem{Def1}[Lemma1]{{Definition}}
\newtheorem{Prop1}[Lemma1]{{Proposition}}
\newtheorem{Claim1}[Lemma1]{{Claim}}
\newtheorem{Rem1}[Lemma1]{{Remark}}
\newtheorem{Cor1}[Lemma1]{{Corollary}}
\newtheorem{Ex1}[Lemma1]{{Example}}
\newtheorem{Not1}[Lemma1]{{Notation}}

\newenvironment{Lemma}{\begin{Lemma1}}{\end{Lemma1}}
\newenvironment{Def}{\begin{Def1}\rm}{\end{Def1}}
\newenvironment{Prop}{\begin{Prop1}}{\end{Prop1}}
\newenvironment{Rem}{\begin{Rem1}\rm}{\end{Rem1}}
\newenvironment{Theorem}{\begin{Theo1}}{\end{Theo1}}

\newenvironment{Cor}{\begin{Cor1}}{\end{Cor1}}

\newenvironment{Example}{\begin{Ex1}\rm}{\end{Ex1}}
\newenvironment{Notation}{\begin{Not1}\rm}{\end{Not1}}

\title[Green correspondence for triangulated categories]{Green correspondence
and relative projectivity for pairs of adjoint functors between triangulated categories}

\author{Alexander Zimmermann}
\address{\newline
Universit\'e de Picardie,
\newline D\'epartement de Math\'ematiques et LAMFA (UMR 7352 du CNRS),
\newline 33 rue St Leu,
\newline F-80039 Amiens Cedex 1,
\newline France}
\email{alexander.zimmermann@u-picardie.fr}

\date{January 28, 2020}

\newcommand{\uar}{\uparrow}
\newcommand{\dar}{\downarrow}
\newcommand{\lra}{\longrightarrow}

\newcommand{\ra}{\rightarrow}
\newcommand{\sdp}{\times\kern-.2em\vrule height1.1ex depth-.05ex}
\newcommand{\epi}{\lra \kern-.8em\ra}

\newcommand{\N}{{\mathbb N}}

\newcommand{\ul}{\underline}
\newcommand{\ol}{\overline}

\newcommand{\Z}{{\mathbb Z}}

\newcommand{\coker}{\textup{coker}}
\newcommand{\add}{\textup{add}}
\newcommand{\thick}{\textup{thick}}
\newcommand{\cone}{\textup{cone}}
\newcommand{\im}{\textup{im}}

\newcommand{\Spec}{\textup{Spec}}
\newcommand{\supp}{\textup{supp}}
\newcommand{\id}{\textup{id}}
\newcommand{\res}{\textup{res}}
\newcommand{\ind}{\textup{ind}}

 \normalbaselines
\subjclass[2010]{Primary: 16E35; Secondary: 18E30;  20C05; 16S34; 18D10 }
\keywords{Green correspondence; Verdier localisation; triangulated category; relative projectivity; vertex of modules}

\begin{document}

\begin{abstract}
Auslander and Kleiner proved in 1994 an abstract version of Green correspondence
for pairs of adjoint functors between three categories. They produce additive quotients
of certain subcategories giving the classical Green correspondence in the special setting of
modular representation theory. Carlson, Peng and Wheeler showed in 1998 that
Green correspondence in the classical setting of modular representation theory is
actually an equivalence between triangulated categories with respect to a non
standard triangulated structure. In the present note we
first define and study a version of relative projectivity, respectively relative injectivity
with respect to pairs of adjoint functors. We then
modify Auslander-Kleiner's construction such that the correspondence holds in the
setting of triangulated categories.
\end{abstract}

\maketitle

\section*{Introduction}

Green correspondence is a very classical and highly important tool in modular representation theory
of finite groups. For a finite group $G$ and a field $k$ of finite characteristic $p$, we associate
to every indecomposable $kG$-module $M$ a $p$-subgroup $D$, called its vertex.
Simplifying slightly, Green correspondence then says that for $H$ being a subgroup of $G$ containing $N_G(D)$,
restriction and induction give a mutually inverse bijection between the indecomposable $kH$-modules
with vertex $D$ and the indecomposable $kG$-modules with vertex $D$. It was known for a long time
that this is actually a categorical correspondence, and in case of trivial intersection Sylow $p$-subgroups it
was known to be more precisely actually an equivalence between the triangulated stable categories.
Only in 1998 Carlson, Peng and Wheeler
showed in \cite{Carlson-Peng-Wheeler} that it is possible to define triangulated structures
also in the general case, and again the Green correspondence is an equivalence between triangulated categories.

Auslander and Kleiner showed in \cite{AuslanderKleiner} that Green correspondence has a
vast generalisation, and actually is a property of pairs of adjoint functors between three categories
$$
\xymatrix{
{\mathcal D}\ar@/^/[r]^{S'}&{\mathcal H}\ar@/^/[l]^{T'}\ar@/^/[r]^{S}&{\mathcal G}\ar@/^/[l]^T
}
$$
such that $(S,T)$ and $(S',T')$ are adjoint pairs and an additional mild hypothesis
on the unit of the adjunction $(S,T)$.
Auslander-Kleiner show that then there is an equivalence between certain additive quotient categories
mimicking the classical Green correspondence. For more details we recall the precise statement as
Theorem~\ref{Greenforadjoints} and Corollary~\ref{Greenfg} below.

Auslander-Kleiner do not study the question whether their abstract Green correspondence
will provide an equivalence between triangulated categories. The present paper aims to fill this gap.
Starting with triangulated categories ${\mathcal D}$, ${\mathcal H}$, ${\mathcal G}$ and pairs of adjoint
triangle functors $(S',T')$ and $(S,T)$ as above, we replace the additive quotient construction
by Verdier localisation modulo the thick subcategories generated by the subcategories for which Auslander and
Kleiner take the additive quotient. We obtain this way triangulated quotient categories and we show the
precise analogue of Theorem~\ref{Greenforadjoints} for the Verdier localisations instead of the additive
quotient categories. In case $S$ is left and right adjoint to $T$, and if in addition the unit of the
adjunction is a monomorphism and the counit is an epimorphism our result shows that
the additive quotient
category is actually already triangulated, and that therefore the
Verdier localisation and the additive quotient
coincide. This way we directly generalise the result of Carlson, Peng and Wheeler
\cite{Carlson-Peng-Wheeler}.

In recent years classification results of thick subcategories of various triangulated categories were
obtained mainly by parameterisations with subvarieties of support varieties. However, most results
use those thick subcategories which form an ideal in an additional monoidal structure, so-called tensor
triangulated categories.
Since many examples, such as non principal blocks of group rings actually are not quite tensor triangulated,
since a unit is missing we study more general a semigroup tensor structure, which is basically
the same as a monoidal structure, but without a unit object. We study properties of our
triangulated Green correspondence in this setting.

We further recall the classical situation and explain how we can recover parts of
the results of Wang-Zhang \cite{Wang-Zhang} and Benson-Wheeler~\cite{Benson-Wheeler} using our approach.

The paper is organised as follows. In Section~\ref{AuslanderKleinerSection}
we recall the main result of Auslander-Kleiner. Generalising the case of relative
projective with respect to subgroups in the case of module categories,
Section~\ref{relativeprojectivetysection}
then introduces the notion of $T$-relative projective (resp. $T$-relatively injective) objects in categories
for functors $T$, and characterises this property in case of $T$ having a left (resp. right) adjoint $S$.
Here we push further and generalise a result due to Brou\'e \cite[Theorem 6.8]{BroueHigman}.
We illustrate our constructions in the case of group algebras.
Section~\ref{VerdierlocalisingSection} then compares Verdier localisation and
the additive quotient categories. We prove there as well our first main result
Theorem~\ref{Greencorrespondencefortriangulated}, generalising
Auslander-Kleiner's theorem to triangulated categories using Verdier localisations.
In Section~\ref{tensortriangulatedsection} we revisit
tensor triangulated categories and study their behaviour within our setting. In particular in
Subsection~\ref{specialcaseofgrouprings} we compare our results to existing results in
the literature in the case of group rings, their stable and derived categories, generalising
various situations in this context.

\subsubsection*{Acknowledgement} I thank the referee for very careful reading and numerous
suggestions which improved greatly the paper, in particular Section 2, but
also throughout. I also thank Olivier Dudas for
giving me the reference  \cite[Theorem 6.8]{BroueHigman}

\section{Summary of Auslander-Kleiner's theory}

\label{AuslanderKleinerSection}

Let ${\mathcal D}$, ${\mathcal H}$, ${\mathcal G}$ be three additive categories and $S,S',T,T'$
be additive functors
$$
\xymatrix{
{\mathcal D}\ar@/^/[r]^{S'}&{\mathcal H}\ar@/^/[l]^{T'}\ar@/^/[r]^{S}&{\mathcal G}\ar@/^/[l]^T
}
$$
such that $(S,T)$ and $(S',T')$ are adjoint pairs. Let $\epsilon:\id_{\mathcal H}\lra TS$
be the unit of the adjunction $(S,T)$. Assume that there is an
endofunctor $U$ of $\mathcal H$ such that $TS=\id_{\mathcal H}\oplus U$, denote by
$p_1:TS\lra \id_{\mathcal H}$
the projection, and suppose that $p_1\circ\epsilon$ is an isomorphism.
If
$\epsilon$ is a split monomorphism, then this is satisfied, but the condition is slightly weaker.
Auslander-Kleiner~\cite{AuslanderKleiner}
prove a Green correspondence result for this situation.

\begin{Notation}\label{inverseimageoffunctors}
\begin{itemize}
\item
For a functor $F:{\mathcal A}\lra{\mathcal B}$ and a full subcategory
$\mathcal V$ of $\mathcal B$ denote for short $F^{-1}({\mathcal V})$ the full
subcategory of ${\mathcal A}$ consisting of objects $A$ such that $F(A)\in add({\mathcal V})$.
\item
For an additive category $\mathcal W$ and an additive subcategory $\mathcal V$ denote by
${\mathcal W}/{\mathcal V}$ the category whose objects are the same objects as those of $\mathcal W$, and
for any two objects $X,Y$ of $\mathcal W$ we put
$$({\mathcal W}/{\mathcal V})(X,Y):={\mathcal W}(X,Y)/I_{\mathcal V}^{\mathcal W}(X,Y),$$
where
$$I_{\mathcal V}^{\mathcal W}(X,Y):=\{f\in{\mathcal W}(X,Y)\;|\;\exists V\in obj({\mathcal V}),
 g\in{\mathcal W}(V,Y),h\in {\mathcal W}(X,V):f=g\circ h\}.$$
\item
If ${\mathcal S}$ and $\mathcal R$ are subcategories of a Krull-Schmidt category $\mathcal W$, then
${\mathcal R}-{\mathcal S}$ denotes the full subcategory of ${\mathcal R}$ consisting of those
objects $X$ of $\mathcal R$ such that no direct factor of $X$ is an object of $\mathcal S$.
\item Recall that a full triangulated subcategory $\mathcal U$ of a triangulated category is
thick (\'epaisse), if it is in addition
closed under taking direct summands (and a fortiori under isomorphisms) in $\mathcal T$.
\item Let $\mathcal U$ be a thick (\'epaisse) subcategory of a triangulated category
$\mathcal T$. Then the Verdier localisation ${\mathcal T}_{\mathcal U}$
(cf \cite[Chapitre I, \S{} 2, no1, no3, no4]{SGA412}, \cite{Verdier})
is the category formed by the same objects as the objects of $\mathcal T$ and
morphisms in  ${\mathcal T}_{\mathcal U}$ are limits of diagrams
$$\xymatrix{X&\ar[l]_-sZ\ar[r]^f&Y}$$  where $f$ and $s$ are morphisms in $\mathcal T$,
and where $\cone(s)$ is an object in $\mathcal U$. The notation we use for the Verdier
localisation is not quite standard,
however in order to distinguish from the additive quotient above we
decided to use this notation (cf Remark~\ref{Remarkdistinguishthenotations}).
\end{itemize}
\end{Notation}

\begin{Theorem} \cite[Theorem 1.10]{AuslanderKleiner} \label{Greenforadjoints}
Assume the hypotheses at the beginning of the section.
Let $\mathcal Y$ be a full additive subcategory of $\mathcal H$ and let
${\mathcal Z}:=(US')^{-1}({\mathcal Y})$.
Then the following two conditions $(\dagger)$ are equivalent.
\begin{itemize}
\item
Each object of $S'T'{\mathcal Y}$ is a direct factor of an object of
$\mathcal Y$ and of an object of $U^{-1}({\mathcal Y})$.
\item Each object of $TSS'T'{\mathcal Y}$ is a direct factor of an object of
$\mathcal Y$.
\end{itemize}
Suppose that the above conditions hold for $\mathcal Y$. Then
\begin{enumerate}
\item $S$ and $T$ induce functors
$${\mathcal H}/S'T'{\mathcal Y}\stackrel{S}{\lra} {\mathcal G}/SS'T'{\mathcal Y}\mbox{ and }
{\mathcal G}/SS'T'{\mathcal Y}\stackrel{T}{\lra}{\mathcal H}/{\mathcal Y}$$
\item For any object $L$ of $\mathcal D$ and any object $B$ of $U^{-1}({\mathcal Y})$ the functor
$S$ induces an isomorphism
$${\mathcal H}/S'T'{\mathcal Y}(S'L,B)\lra {\mathcal G}/SS'T'{\mathcal Y}(SS'L,SB)$$
\item
For any object $L$ of $(US')^{-1}{\mathcal Y}$ and any object $A$ of ${\mathcal G}$ the functor
$T$ induces an isomorphism
$${\mathcal G}/SS'T'{\mathcal Y}(SS'L,B)\lra {\mathcal H}/{\mathcal Y}(TSS'L,TA)$$
\item
The restrictions of $S$
$$(\add S'{\mathcal Z})/S'T'{\mathcal Y}\stackrel S\lra (\add SS'{\mathcal Z})/SS'T'{\mathcal Y}$$
and $T$
$$(\add SS'{\mathcal Z})/SS'T'{\mathcal Y}\stackrel T\lra (\add TSS'{\mathcal Z})/{\mathcal Y}$$
are equivalences of categories, and
$$(\add S'{\mathcal Z})/S'T'{\mathcal Y}\stackrel{TS}{\lra}(\add TSS'{\mathcal Z})/{\mathcal Y}$$
is isomorphic to the natural projection.
\item
If each object of $S'T'US'{\mathcal D}$ is a direct factor of $US'{\mathcal D}$, then
${\mathcal Y}=US'{\mathcal D}$ satisfies the hypothesis of the theorem.
\end{enumerate}
\end{Theorem}

A main consequence is

\begin{Cor} \cite[Corollary 1.12]{AuslanderKleiner}\label{Greenfg}
Let $\mathcal Y$ be a full additive subcategory of $\mathcal H$ satisfying $(\dagger)$ of
Theorem~\ref{Greenforadjoints}
and suppose that $\mathcal H$ and $\mathcal G$ are both Krull-Schmidt categories.
Using the notations of Theorem~\ref{Greenforadjoints} then the following hold.
\begin{enumerate}
\item
For each indecomposable object $N$ of $(add(S'{\mathcal Z}))-{S'T'{\mathcal Y}}$ the object $SN$
has a unique
indecomposable direct factor $g(N)$ which is not a direct factor of an object in $SS'T'{\mathcal Y}$.
\item
For each indecomposable object $M$ of $(add(SS'{\mathcal Z}))-{SS'T'{\mathcal Y}}$
the object $TM$ has a unique
indecomposable direct factor $f(M)$ which is not a direct factor of an object in ${\mathcal Y}$.
\item
$f(g(N))=N$
\item
$g(f(M))=M$.
\end{enumerate}
\end{Cor}

\section{Relative projectivity and injectivity with respect to pairs of adjoint functors}

\label{relativeprojectivetysection}

\subsection{Relative homological algebra revisited}

\label{relativehomologicalalgebrarevisited}

We shall need to revise some facts from relative homological
algebra, following \cite{BeligiannisMarmaridis}.
Recall that a full subcategory $\mathcal X$ of an additive category
$\mathcal S$ is contravariantly finite if
for any object $S$ of $\mathcal S$ there is an object $X$ of $\add \mathcal X$ and a morphism
$f\in{\mathcal S}(X,S)$ such that for any $X'$ in $\mathcal X$ the induced map
$${\mathcal S}(X',f):{\mathcal S}(X',X)\lra {\mathcal S}(X',S)$$
is surjective. We call such an object $X$ of $\mathcal S$ a right ${\mathcal X}$-approximation of $S$.
The dual notion, using the covariant $Hom$-functor leads to the notion
of a covariantly finite subcategory.
With this notion in mind we shall have the following

\begin{Lemma} (Auslander-Reiten \cite[Proposition 1.2]{AR})
If the additive functor $T:{\mathcal S}\rightarrow{\mathcal T}$
between the additive categories ${\mathcal S}$ and ${\mathcal T}$
admits a left adjoint $S_\ell$, then
$\add(\im(S_\ell))$ is a contravariantly finite subcategory of $\mathcal S$.
If $T$ admits a right adjoint $S_r$, then $\add(\im(S_r))$ is a
covariantly finite subcategory of $\mathcal S$.
\end{Lemma}

\begin{proof}
Let ${\mathcal X}:=\add(\im(S_\ell))$.
Consider the counit
$$\eta:S_\ell T\lra \id_{\mathcal S}$$
of the adjoint pair $(S_\ell,T)$. Evaluation on any object $Q$ of $\mathcal S$
gives a morphism
$$\eta_Q:S_\ell T Q\lra Q.$$
Now, given an object $S_\ell P$ in $\im(S_\ell)$, we have
$$\xymatrix{{\mathcal S}(S_\ell P,S_\ell T Q)\ar[rr]^-{{\mathcal S}(S_\ell P,\eta_Q)}&&
{\mathcal S}(S_\ell P, Q)\\
{\mathcal T}( P,TS_\ell T Q)\ar[u]_{\simeq}&&{\mathcal T}( P, T Q)\ar[u]^{\simeq}}$$
and the composite map is a split epimorphism by \cite[IV Theorem 1.(ii)]{Maclane}.
Since the property holds true for direct factors of an object $S_\ell P$ as well,
we showed that ${\mathcal X}$ is a contravariantly
finite subcategory of $\mathcal S$.

By the dual argument, if $T$ admits a right adjoint $S_r$, then $\add(\im S_r)$ is
a covariantly subcategory of $\mathcal S$.
\end{proof}

Recall from Beligiannis and Marmaridis~\cite{BeligiannisMarmaridis} that we may produce from
contravariantly finite subcategories a relative homological algebra.
Let $\mathcal X$ be a contravariantly finite subcategory of an additive category $\mathcal S$.
Then a morphism $g\in{\mathcal S}(A,B)$ is said to be  ${\mathcal X}$-epic
if for any object $X$ of ${\mathcal X}$ the morphism
$${\mathcal S}(X,g):{\mathcal S}(X,A)\lra{\mathcal S}(X,B)$$
is surjective. By the very definition, a right ${\mathcal X}$-approximation
is an ${\mathcal X}$-epic. If $\mathcal X$ is contravariantly finite and if each
$\mathcal X$-epic
has a kernel, \cite[Theorem 2.12]{BeligiannisMarmaridis} shows that for any object
$S$ of $\mathcal S$ the choice of a right $\mathcal X$-approximation $X_S\ra S$
induces a left triangulation on the stable category ${\mathcal S}/{\mathcal X}$. Moreover, two
such choices give equivalent left triangulated categories. Hence, a contravariantly finite
subcategory $\mathcal X$ such that each $\mathcal X$-epic
has a kernel gives rise to the relative $Ext^n$-group with respect to $\mathcal X$,
denoted by $Ext^n_{\mathcal X}(A,B)$ namely the
evaluation on the object $B$, of the $n$-th derived functor of ${\mathcal S}(-,B)$,
obtained by a $\mathcal X$-resolution of $A$.

\begin{Lemma}\label{imepicsplitepi}
Let $\mathcal S$ and $\mathcal T$ be additive categories, let $T:{\mathcal S}\ra{\mathcal T}$
be an additive functor admitting a left adjoint $S_\ell$. Then
$g\in{\mathcal S}(A,B)$ is $\add(\im(S_\ell))$-epic if and only if $T(g)$ is a split epimorphism.
\end{Lemma}

\begin{proof}
Let $g\in{\mathcal S}(A,B)$ be $\add(\im(S_\ell))$-epic. Then for any object $C$ of
$\mathcal T$ we get
$${\mathcal S}(S_\ell(C),g):{\mathcal S}(S_\ell(C),A)\lra {\mathcal S}(S_\ell(C),B)$$
is surjective. Hence, for any object $C$ of $\mathcal T$ we get
$${\mathcal T}(C,Tg):{\mathcal T}(C,TA)\lra {\mathcal T}(C,TB)$$
is surjective. In particular, for $C=TB$ there is $f\in {\mathcal T}(TB,TA)$
with $Tg\circ f=\id_{TB}$. Hence $Tg$ is a split epimorphism.

Let $g\in{\mathcal S}(A,B)$ be such that $Tg$ is a split epimorphism. Then there is
$f\in {\mathcal T}(TB,TA)$ with $Tg\circ f=\id_{TB}$. Let $C$ be an object of $\mathcal T$ and
let $h\in {\mathcal S}(S_\ell(C),B)$.
We need to show that there is $k\in{\mathcal S}(S_\ell C,A)$
such that ${\mathcal S}(S_\ell(C),g)(k)=g\circ k=h$, where
$${\mathcal S}(S_\ell(C),g):{\mathcal S}(S_\ell(C),A)\lra {\mathcal S}(S_\ell(C),B).$$
Since $T$ is right adjoint to $S_\ell$, this is equivalent to
$${\mathcal T}(C,Tg):{\mathcal T}(C,TA)\lra {\mathcal T}(C,TB)$$
is surjective. For $h\in{\mathcal T}(C,TB)$ we get
$$h=(Tg\circ f)\circ h=Tg\circ (f\circ h)={\mathcal T}(C,Tg)(f\circ h)={\mathcal S}(S_\ell C,g)(f\circ h).$$
Clearly we can pass to direct factors of $S_\ell C$. Therefore, $g$ is $\add(\im( S_\ell))$-epic.
\end{proof}

Note that Lemma~\ref{imepicsplitepi} has a dual version for functors $T$ admitting right
adjoint functors $S_r$.

Again, in the setting of \cite{BeligiannisMarmaridis} the
$\add(\im(S_\ell))$-relative projectives, are those objects $Q$
with $$Ext^n_{\add(\im S_\ell)}(Q,B)=0$$ for all objects $B$ and $n>0$. By definition, this
coincides with the objects $Q$ for which the counit of the adjunction $(S_\ell, T)$ splits.
These are precisely the objects in $\add(\im(S_\ell))$.

The dual statement applies in case of $T$ having a right adjoint $S_r$, and
considering covariantly finite subcategories and $\add(\im(S_r))$-coresolutions
instead of contravariantly finite subcategories and $\add(\im(S_\ell))$-resolutions.

This motivates the following definition.

\begin{Def}\label{relativeprojectivitydefinition}
Let ${\mathcal T}$ and ${\mathcal S}$ be triangulated categories, and let
$T:{\mathcal S}\ra{\mathcal T}$
be a triangle functor.
\begin{itemize}
\item
Suppose $T$ has a left adjoint.
Then an object $Q$ of $\mathcal S$ is
{\em  $T$-relative projective} if the natural
transformation
$${\mathcal S}(Q,-)\ra {\mathcal T}(TQ,T-)$$
induced by $T$ is injective.
\item
Suppose $T$ has a right adjoint.
Then an object $Q$ of $\mathcal S$ is
{\em  $T$-relative injective} if the natural
transformation
$${\mathcal S}(-,Q)\ra {\mathcal T}(T-,TQ)$$
induced by $T$ is injective.
\end{itemize}
\end{Def}

\begin{Rem}
Recall that for a field $k$ of finite characteristic $p>0$ and a finite group $G$
with a subgroup $H$, an indecomposable
$kG$-module $M$ is called relatively $H$-projective if each epimorphism
$N\twoheadrightarrow M$ of $kG$-modules, which is known to be
split as $kH$-module morphism, splits as $kG$-module morphism.
This definition of relative projectivity was developed by Hochschild \cite{Hochschild}
in the situation of a ring $R$, a subring $S$ of $R$. Hochschild declares
an $R$-module $M$ to be $(S,R)$projective, if any short exact sequence
$$0\lra X\lra Y\lra M\lra 0$$ which is known to be split as short exact sequence of
$S$-modules, is automatically split as short exact sequence of $R$-modules.
Denoting by $\res^R_S:R-Mod\lra S-Mod$ the exact functor given by restriction to the subring $S$,
this translates into slightly more modern terms into the statement that
$M$ is $(S,R)$-projective if and only if
$$Ext^1_R(M,X)\lra Ext^1_S(\res^R_S(M),\res^R_S(X))$$
is injective for any $R$-module $X$.
Hence, since $\res^R_S(X)[1]\simeq \res^R_S(X[1])$ we get that
$M$ is $(S,R)$-projective if and only if
$$Hom_{D^b(R)}(M,X[1])\lra Hom_{D^b(S)}(\res^R_S(M),\res^R_S(X[1]))$$
is injective for all objects $X$. Since each object $X$ can be seen as an object $X=Y[-1]$,
Definition~\ref{relativeprojectivitydefinition} could make sense in a broader
context. We will not elaborate on this here (cf \cite{RemarksOnGreen}).
\end{Rem}

\begin{Rem}
Grime \cite{Grime} defines an object to be
relative projective with respect to a functor $F$ admitting a left adjoint $L$
as those which are direct factors of an object in the image of $L$.
This is a direct generalisation of Green's definition \cite{Green}, whereas
our definition is closer to Hochschild's definition \cite{Hochschild}.
However,
the concepts coincide, as will be shown in
Proposition~\ref{HigmansLemmaForAdjointFunctors} below.
\end{Rem}

\subsection{Relative projectivity for triangulated categories}

\label{relativeprojectivityfortriangulatedsection}

We shall study the concept of $T$-relative projectivity/injectivity from
Definition~\ref{relativeprojectivitydefinition} for
triangle functors $T$ between triangulated categories
admitting a left adjoint $S_\ell$
and a right adjoint $S_r$. Then the concept has a very nice interpretation.

\begin{Lemma}\label{relatively-t-projective-for-adjoints}
Let ${\mathcal S}$ and ${\mathcal T}$ be additive categories and let
$T:{\mathcal S}\ra{\mathcal T}$ be an additive functor.
\begin{itemize}
\item
If $T$ has a left adjoint $S_\ell$, then an object $Q$ is $T$-relative projective
if and only if the evaluation on $Q$ of
the counit $\eta$ of the adjuntion  $\eta_Q:S_\ell TQ\ra Q$ is an epimorphism.
Any object in $\add(\im(S_\ell))$ is $T$-relative projective.
\item
If $T$ has a right adjoint $S_r$, then an object $Q$ is $T$-relative injective
if and only if the evaluation on $Q$ of
the unit $\eta$ of the adjuntion $\epsilon_Q:Q\ra S_rTQ$ is a monomorphism.
Any object in $\add(\im(S_r))$ is $T$-relative injective.
\end{itemize}
\end{Lemma}

\begin{proof}
Suppose that $T$ has a left adjoint $S_\ell$.
Then the conuit $\eta_Q:S_\ell TQ\ra Q$ is an epimorphism if and only if for any object $A$
the morphism $${\mathcal S}(\eta_Q,A):{\mathcal S}(Q,A)\ra {\mathcal S}(S_\ell TQ,A)$$
is a monomorphism. This in turn is equivalent to the statement that the natural
transformation of functors ${\mathcal S}\ra \Z-Mod$
$${\mathcal S}(\eta_Q,-):{\mathcal S}(Q,-)\ra {\mathcal S}(S_\ell TQ,-)$$
is a monomorphism. Using the defining property of $(S_\ell, T)$ being an adjoint pair,
this is equivalent to
$${\mathcal S}(\eta_Q,-):{\mathcal S}(Q,-)\ra {\mathcal T}(TQ,T-)$$
being a monomorphism. Hence, the statement is equivalent to $Q$ being $T$-relative projective.
Now $\eta_{S_\ell Q'}$ is a split epimorphism  for any object $Q'$ of
$\mathcal T$ by \cite[IV Theorem 1.(ii)]{Maclane}.

Suppose that $T$ has a right adjoint $S_r$.
Then the unit $\epsilon_Q:Q\ra S_rTQ$ is a monomorphism if and only if
$${\mathcal S}(A,\epsilon_Q):{\mathcal S}(A,Q)\ra {\mathcal S}(A,S_r TQ)$$
is a monomorphism. This is equivalent to
$${\mathcal S}(-,\epsilon_Q):{\mathcal S}(-,Q)\ra {\mathcal S}(-,S_r TQ)={\mathcal T}(T-,TQ)$$
is a monomorphism, which is equivalent to $Q$ is $T$-relative injective.
Now $\epsilon_{S_r Q'}$ is a split monomorphism for any object $Q'$ of
$\mathcal T$ by \cite[IV Theorem 1.(ii)]{Maclane}.
\end{proof}

\begin{Rem}\label{relatively-t-projective-triangulated}
Note that in a triangulated category $\mathcal S$ the notions of epimorphism
(resp. monomorphism) and split epimorphism (resp. split monomorphism) coincide.
\end{Rem}

\begin{Prop}\label{relativelyprojectivecharacterisation}
Let ${\mathcal T}$ and ${\mathcal S}$ be triangulated categories and let
$T:{\mathcal S}\lra {\mathcal T}$
be a triangle functor. Suppose that $T$ has a left (respectively right) adjoint $S$.
Then an object $Q$ is $T$-relative projective (respectively injective)
if and only if $Q$ is in $\add(\im(S))$.
\end{Prop}

\begin{proof}
By Lemma~\ref{relatively-t-projective-for-adjoints} and Remark~\ref{relatively-t-projective-triangulated}
$Q$ is $T$-relative projective (respectively $T$-relative injective)
if and only if $Q$ is in $\add(\im S)$.
\end{proof}

\begin{Cor}\label{counitsplitsmeansrelativelyprojective}
Let ${\mathcal T}$ and ${\mathcal S}$ be triangulated categories and let
$T:{\mathcal S}\lra {\mathcal T}$
be a triangle functor. Suppose that $T$ has a left (respectively right)
adjoint $S$, and let $\eta:ST\lra\textup{id}$ be the
counit (respectively $\widetilde\epsilon:\textup{id}\lra ST$ the unit) of the adjunction.
Then $Q$ is $T$-relative projective (respectively injective) if and
only if $\eta_Q$ is a split epimorphism
(respectively $\widetilde\epsilon_Q$ is a split monomorphism).
\end{Cor}

\begin{proof}
This is precisely Proposition~\ref{relativelyprojectivecharacterisation} in connection with
Lemma~\ref{relatively-t-projective-for-adjoints} and Remark~\ref{relatively-t-projective-triangulated}.
\end{proof}

We summarize the situation to an analogue of Higman's lemma for pairs of adjoint functors
between triangulated categories.

\begin{Prop}\label{HigmansLemmaForAdjointFunctors}
Let ${\mathcal T}$ and ${\mathcal S}$ be triangulated categories and let $T:{\mathcal S}\lra {\mathcal T}$
be a triangle functor. Suppose that $T$ has a left (respectively right) adjoint $S$. Let $M$ be an
 object of $\mathcal T$. Then the following
are equivalent:
\begin{enumerate}
\item $M$ is $T$-relative projective (respectively injective).
\item $M$ is in $\add(\im S)$.
\item $M$ is a direct factor of some $S(L)$ for some $L$ in $\mathcal S$.
\item $M$ is a direct factor of some $ST(M)$.
\end{enumerate}
\end{Prop}

\begin{proof}
$(1)\Leftrightarrow (2)$ is Proposition~\ref{relativelyprojectivecharacterisation}.

$(2)\Leftrightarrow (3)$ is the definition of $\add(\im S)$.

$(3)\Rightarrow (4)$ is trivial.

$(4)\Rightarrow (1)$ is Corollary~\ref{counitsplitsmeansrelativelyprojective}.
\end{proof}

\begin{Rem}
Note that Corollary~\ref{counitsplitsmeansrelativelyprojective} generalises
\cite[Proposition 2.1.6, Proposition 2.1.8]{reptheobuch} to this more general situation.
\end{Rem}

\begin{Rem}
In case $T$ has a left adjoint $S$, which is also assumed to be a right adjoint,
and $\mathcal T$ is an abelian or  triangulated category
Brou\'e defined $T$-relative projective and $T$-relative injective objects
in \cite[Theorem 6.8]{BroueHigman}. In this situation he showed a version of Higman's lemma
as Proposition~\ref{HigmansLemmaForAdjointFunctors} by completely different means.
\end{Rem}

\begin{Cor}
Let $\mathcal S$ and $\mathcal T$ be triangulated categories and
let $T:{\mathcal S}\lra{\mathcal T}$ be a triangle functor admitting
a left (resp. right) adjoint $S$. Then all objects of $\mathcal S$ are
$T$-relative projective (resp. injective) if and only if
${\mathcal S}=\add(\im(S))$.
\end{Cor}

Note that all objects of
$\mathcal S$ are
$T$-relative projective (resp. injective) if and only if the global
dimension of the
relative homological algebra described in Section~\ref{relativehomologicalalgebrarevisited}
is $0$.

\begin{Example}\label{examplederivedcategoriesofgroups}
Let $G$ be a group, and let $H$ be a subgroup of finite index $n$.
Denote by $\dar^G_H$ the functor given by restriction of the $G$-action
to the $H$-action, and by $\uar_H^G$ the functor $kG\otimes_{kH}-$ given by induction.
If $n$ is invertible in the field $k$,
then every object $M$ in $D^b(kG)$ is $\dar^G_H$-relative projective.
Indeed, \cite[Proposition 2.1.10]{reptheobuch} shows that the multiplication
$kG\otimes_{kH}kG\lra kG$ splits as morphism of $kG-kG$-bimodules.
The counit of the adjunction
$(\uar_H^G,\dar_H^G)$ is $kG\otimes_{kH}kG\otimes_{kG}-\lra kG\otimes_{kG}-$,
and by hypothesis
this map splits.
\end{Example}

\subsection{Revisiting the case studied by Carlson-Peng-Wheeler}

The purpose of this section is to give a structural explanation of an argument
in the proof of Carlson, Peng and Wheeler for the statement that
the relative stable category is triangulated
(cf \cite[page 304; proof of Theorem 6.2]{Carlson-Peng-Wheeler}).
Note that Grime gave a slightly less general
structural explanation in \cite[Example 3.6]{Grime}.

\begin{Rem}\label{whendowehavetriangulatedadditivequotients}
Carlson, Peng and Wheeler consider the classical case of group rings,
namely let $k$ be a field of characteristic $p>0$, let $G$ be a finite group,
let $D$ be a $p$-subgroup of $G$ and let $H$ be a subgroup of $G$ containing the normalizer of $D$ in $G$.
They consider
the additive quotient of the module category modulo the
morphisms which factor through $\dar^G_E$-projective modules, for some $E\in\mathfrak Y$, where
$\mathfrak Y=\{E\leq H\cap D^g\;|\;g\in G\setminus H\}$, and show that this produces a
triangulated category. Carlson, Peng and Wheeler
use the general approach given by Happel \cite[Theorem I.2.6]{Happel}
which shows that the additive quotient of any Frobenius category modulo
relative injective-projectives is triangulated.
However, Carlson, Peng and Wheeler just mention that Happel's proof for
Frobenius categories to have triangulated stable categories carries over to
this more general situation. The purpose of this section is to
show that the fact that the proof carries over has a structural reason, and uses
more precisely the properties from Section~\ref{relativehomologicalalgebrarevisited}
and Section~\ref{relativeprojectivityfortriangulatedsection}.
\end{Rem}

Note that group rings are symmetric, hence the module category is Frobenius. Moreover,
the functors considered in classical Green correspondence, namely restriction and induction,
are left and right adjoint to each other.
We note that in our general abstract situation relative injectives and relative projectives do not
coincide in general. The situation changes in case $S$
is at the same time left and right adjoint to $T$
and the categories are already Frobenius categories.

\medskip

Recall from \cite[Definition 2.1]{Buehler} the concept of an exact category.
Let $\mathcal A$ be an additive category. Given three objects $A_1,A_2,A_3$ in
$\mathcal A$ and  $f\in{\mathcal A}(A_1,A_2)$ and $g\in {\mathcal A}(A_2,A_3)$.
Then $(f,g)$ is a short exact sequence, denoted by
$$\xymatrix{0\ar[r]& A_1\ar[r]^{f}&A_2\ar[r]^g& A_3\ar[r]& 0,}$$
(or occasionally by $\xymatrix{ A_1\ar@{^{(}->}[r]^{f}&A_2\ar@{->>}[r]^g& A_3,}$)
if $\ker(g)=f$ and $g=\coker(f)$.

An exact structure on the additive category $\mathcal A$ is given by a class
$E_{\mathcal A}$ of short exact sequences, called admissible short exact sequences,
satisfying the following axioms below. If
$$\xymatrix{0\ar[r]& A_1\ar[r]^{f}&A_2\ar[r]^g& A_3\ar[r]& 0,}$$
is a short exact sequence in $E_{\mathcal A}$, then we say that $f$ is an
admissible monomorphism and $g$ is an admissible epimorphism.
\begin{itemize}
\item For all objects $A$ the identity on $A$ is admissible monomorphism and admissible epimorphism.
\item Admissible monomorphisms are closed under composition, and
admissible epimorphisms are closed under composition.
\item If $\alpha:X\ra Y$ is an admissible monomorphism, and $f:X\ra Z$ is any morphism, then
the pushout $$\xymatrix{X\ar[r]^\alpha\ar[d]_f&Y\ar[d]^{\check f}\\ Z\ar[r]^{\check\alpha}&U}$$
exists and $\check\alpha$ is an admissible monomorphism.
\item If $\alpha:Y\ra X$ is an admissible epimorphism, and $f:Z\ra X$ is any morphism, then
the pullback $$\xymatrix{Y\ar[r]^\alpha&X\\ U\ar[r]^{\hat\alpha}\ar[u]^{\hat f}&Z\ar[u]_f}$$
exists and $\hat\alpha$ is an admissible epimorphism.
\end{itemize}
An exact category is an additive
category $\mathcal A$ with a class $E_{\mathcal A}$ of short exact sequences,
stable under isomorphism and satisfying the above axioms.
See \cite{Buehler} for an exhaustive development of exact categories.

\begin{Prop}\label{frobeniusstructuresforfunctors}
Let $({\mathcal S},E_{\mathcal S})$ and $({\mathcal T},E_{\mathcal T})$ be exact categories
with $E_{\mathcal S}$ respectively $E_{\mathcal T}$ being the class of admissible exact sequences.
If $T:{\mathcal S}\lra{\mathcal T}$ is a  functor with a left adjoint
$S_\ell$ and a right adjoint $S_r$,
\begin{itemize}
\item
then
$$E_T:=\left\{\left(\xymatrix{X\ar@{^{(}->}[r]^f& Y\ar@{->>}[r]^g& Z}\right)\in E_{\mathcal S}\;|
\;\left(\xymatrix{TX\ar@{^{(}->}[r]^{Tf}&TY\ar@{->>}[r]^{Tg}&TZ}\right)\in E_{\mathcal T}\right\}$$
defines an exact structure on $\mathcal S$.
\item
If moreover the unit $\epsilon:\id_{\mathcal S}\lra S_rT$ is an admissible monomorphism in
$E_{\mathcal S}$
and if the counit $\eta:S_\ell T\lra \id_{\mathcal S}$ is an admissible epimorphism in
$E_{\mathcal S}$,
\begin{itemize}
\item
then $({\mathcal S},E_T)$ has enough $T$-relative projectives and enough
$T$-relative injectives.
\item Suppose now in addition that $\mathcal S$ and $\mathcal T$ are abelian
Frobenius (i.e. an abelian category which is Frobenius with respect
to the class of all exact sequences).
Then the class of
$T$-relative projectives coincides with the class of objects in $\add(\im S_\ell)$
and the class of $T$-relative injectives coincides with the class of
objects in $\add(\im S_r)$.
\end{itemize}
\end{itemize}
\end{Prop}

\begin{proof}
We first show that $E_T$ is an exact structure.
$T(\id_A)=\id_{TA}$, which implies the first condition.
$T$ maps compositions to compositions, and hence compositions of
admissible monics/epics are admissible monics/epics. Sequences are closed under
isomorphisms, as $T$ is exact and hence maps isomorphisms to isomorphisms.
Let $\xymatrix{X\ar[r]^\alpha&Y\ar[r]^\beta&Z}$
be an exact sequence in $E_{\mathcal S}$ and let $\xymatrix{X\ar[r]^{f}&X'}$ be any morphism.
Then, since $E_{\mathcal S}$ is an exact structure, we may form the pushout
$$\xymatrix{X\ar[r]^\alpha\ar[d]^{f}&Y\ar[r]^\beta\ar[d]^g&Z\ar@{=}[d]\\
X'\ar[r]^{\check\alpha}&Y'\ar[r]^{\check \beta}&Z}.$$
As $E_{\mathcal S}$ is an exact structure, the lower row is an element of $E_{\mathcal S}$.
The sequence $$\xymatrix{0\ar[r]&X\ar[r]^-{f\choose {-\alpha}}&
X'\oplus Y\ar[r]^-{(\check\alpha\;,\;g)}&Y'\ar[r]&0}$$
is exact, since the above is a pushout and $\alpha,\check\alpha$ are monomorphisms.
Since $T$ is exact,
$$\xymatrix{0\ar[r]&TX\ar[r]^-{Tf\choose {-T\alpha}}&
TX'\oplus TY\ar[rr]^-{(T\check\alpha\;,\;Tg)}&&TY'\ar[r]&0}$$
is exact. Therefore
$$\xymatrix{TX\ar[r]^-{T\alpha}\ar[d]^{Tf}&TY\ar[r]^-{T\beta}\ar[d]^{Tg}&TZ\ar@{=}[d]\\
TX'\ar[r]^{T\check\alpha}&TY'\ar[r]^{T\check \beta}&TZ}$$
is a pushout diagram. Since the above row is in $E_{\mathcal T}$, and since $E_{\mathcal T}$
is an exact structure, also the lower row is in $E_{\mathcal T}$.
This shows the third axiom. Dually also the
fourth axiom holds.

\medskip

We now assume the additional condition on the unit and the counit.
The fact that $\add(\im S_\ell)$ are $T$-relative injective objects and $\add(\im S_r)$
are $T$-relative projective
objects is Lemma~\ref{relatively-t-projective-for-adjoints}.
The fact that we then get enough $T$-relative projective objects follows from the
hypothesis on
the counit, and the fact that we then get enough $T$-relative injective objects
follows from the hypothesis on
the unit.

\medskip

The hypothesis on $\mathcal S$ and $\mathcal T$ being Frobenius with respect to all
short exact sequences implies that the stable categories
modulo projective-injective objects $\ul{\mathcal S}$ and $\ul{\mathcal T}$ are
triangulated (cf Happel~\cite[Theorem I.2.6]{Happel}).
Proposition~\ref{relativelyprojectivecharacterisation} applied to this triangulated category
shows that $\add(\im S_\ell)$ are precisely the $T$-relative projective objects
and $\add(\im S_r)$ are precisely the $T$-relative injective objects of this
new exact structure.
\end{proof}

\begin{Rem}
Note that the hypothesis of $\epsilon$ being a monomorphism and $\eta$ being an epimorphism
for the adjunctions involved is very strong.
For an abelian category $\mathcal S$,
if $T:{\mathcal S}\ra{\mathcal T}$ has left and right adjoints $S_\ell$ and $S_r$, then $T$ is exact.
Further, Eilenberg-Moore~\cite[Proposition 1.5]{Eilenberg-Moore}
(cf also Grime \cite[Lemma 2.1]{Grime}) show that the
unit $\id\ra S_r T$, as in Proposition~\ref{frobeniusstructuresforfunctors}, is a
monomorphism if and only if $TX=0$ implies $X=0$,
if and only if the counit $S_\ell T\ra \id$ is an epimorphism.
The counit $S_\ell T\ra \id$ is an epimorphism
if and only if $T$ is faithful.

In order to assure all quotient categories in Theorem~\ref{Greenforadjoints} being triangulated,
using Proposition~\ref{frobeniusstructuresforfunctors} we need to
assume the hypothesis for all the functors $S,S',T,T'$, and hence get quite a few restrictions on
these functors.
\end{Rem}

\begin{Rem}
The first item in Proposition~\ref{frobeniusstructuresforfunctors}
should be compared with the statement \cite[Theorem II.2.1]{Eilenberg-Moore} by Eilenberg-Moore.
\end{Rem}

\begin{Rem}\label{Frobeniusforleftequalsrightadjoint} Let
$({\mathcal S},E_{\mathcal S})$ and $({\mathcal T},E_{\mathcal T})$ be exact categories,
let $T:{\mathcal S}\ra{\mathcal T}$ be an exact functor admitting a left
adjoint $S_\ell$ and a right
adjoint $S_r$, and suppose the unit $\epsilon:\id_{\mathcal S}\lra S_rT$ of the adjoint property
$(T,S_r)$
is a monomorphism, and the counit $\eta:S_\ell T\lra \id_{\mathcal S}$ of the adjoint
property $(S_\ell,T)$
is an epimorphism. Suppose moreover that $\add(\im(S_r))=\add(\im(S_\ell))$.
If in addition $\mathcal S$ and $\mathcal T$ are abelian Frobenius
categories (i.e. an abelian category which is Frobenius
with respect to the exact structure given by all exact sequences,
then Proposition~\ref{frobeniusstructuresforfunctors}
shows that
$$E_T:=\left\{\left(\xymatrix{X\ar@{^{(}->}[r]^f& Y\ar@{->>}[r]^g& Z}\right)\in E_{\mathcal S}\;|
\;\left(\xymatrix{TX\ar@{^{(}->}[r]^{Tf}&TY\ar@{->>}[r]^{Tg}&TZ}\right)\in E_{\mathcal T}\right\}$$
is a Frobenius structure on $\mathcal S$. We call this the $T$-relative Frobenius structure.
Following Happel~\cite[Theorem I.2.6]{Happel} the stable category $\ul{\mathcal S}^T$ of
$\mathcal S$ modulo the $T$-relative projectives is in this case
a triangulated category. The distinguished triangles are
constructed as follows.
Given two objects $M$ and $N$ in $\mathcal S$ and  $f\in{\mathcal S}(M,N)$.
Then we may form the pushout diagram
$$
\xymatrix{
0\ar[r]&M\ar[r]^{\widetilde\epsilon_M}\ar[d]_f&S_rTM\ar[r]\ar[d]&\Omega_T^{-1}(M)\ar[r]\ar@{=}[d]&0\\
0\ar[r]&N\ar[r]^{c_1(f)}&C(f)\ar[r]^{{c_2(f)}}&\Omega_T^{-1}(M)\ar[r]&0
}
$$
(or analogously the pullback diagram along $\eta_N:S_\ell TN\epi N$).
Then $\underline{\mathcal S}^T$ is a triangulated category with distinguished triangles being
isomorphic to triangles of the form
$$\xymatrix{M\ar[r]^f&N\ar[r]^{c_1(f)}&C(f)\ar[r]^{c_2(f)}&\Omega_T^{-1}(M)}$$
for any $f\in{\mathcal S}(M,N)$.
\end{Rem}

We recall a result implicit in Grime \cite{Grime}.

\begin{Prop} (cf \cite[Theorem 3.3]{Grime})
Let $({\mathcal S},E_{\mathcal S})$ and $({\mathcal T},E_{\mathcal T})$
be exact categories and let $T:{\mathcal S}\ra{\mathcal T}$
be a functor which admits a left adjoint $S_\ell$ and a right adjoint $S_r$.
Assume that the counit $S_\ell T\ra \id_{\mathcal S}$ of the adjoint pair $(S_\ell, T)$
is an admissible epimorphism in the exact category $({\mathcal S},E_{\mathcal S})$ and that
the unit $\id_{\mathcal S}\ra S_rT$ of the adjoint pair $(T,S_r)$ is an admissible monomorphism
in the exact category $({\mathcal S},E_{\mathcal S})$.
Putting
$$E_T:=\{\left(\xymatrix{ X\ar@{^{(}->}[r]^f& Y\ar@{->>}[r]^{g}& Z}\right)\in E_{\mathcal S}\;|\;
\xymatrix{TX\ar@{^{(}->}[r]^{Tf}& TY\ar@{->>}[r]^{Tg}& TZ}\;\textup{ is split exact in }{\mathcal T}\}$$
then $({\mathcal S},E_T)$ is an exact category with enough projective and enough injective
objects. The full subcategory of projective objects coincides with the
full subcategory $\add(\im(S_\ell))$ and the full subcategory of injective
objects coincides with the full subcategory $\add(\im(S_r))$.
\end{Prop}

Grime's proposition follows from Proposition~\ref{frobeniusstructuresforfunctors}
when it is applied to the case of the split exact structure on $\mathcal T$.

\subsection{Relative projectivity for derived categories of group rings}

We shall apply our concept of relative projectivity to the special case of
the derived category of a block of a group ring $kG$. We first note that
if $A$ is a  finite dimensional $k$-algebra over a field $k$, then
$D^b(A)$ is a Krull-Schmidt category.
Let $G$ be a finite group, let $H$ be a subgroup of $G$,
let $k$ be a field of characteristic $p>0$.
Then we consider the functors $\uar_H^G$ and $\dar^G_H$. Note that
both functors are exact functors between $kG-mod$ and $kH-mod$.
These functors form an adjoint pair, in the sense that $(\uar_H^G,\dar^G_H)$ and
$(\dar_H^G,\uar^G_H)$ are both adjoint pairs. Note that since
both functors are exact, they provide functors
$$S:=\uar_H^G:D^b(kH)\lra D^b(kG)$$
and
$$T:=\dar_H^G:D^b(kG)\lra D^b(kH).$$
We define ${\mathcal G}:=D^b(kG)$ and ${\mathcal H}:=D^b(kH)$.
Moreover, $(\uar_H^G,\dar^G_H)$ and $(\dar_H^G,\uar^G_H)$ are both adjoint pairs
also between the derived categories. As for its restriction to the module categories
we have

\begin{Lemma}\label{induceistransitive}
Let $K\leq H\leq G$ be an increasing sequence of groups.
Then for the functors
\begin{eqnarray*}
\uar_H^G&:&D^b(kH)\lra D^b(kG),\\
\uar_K^H&:&D^b(kK)\lra D^b(kH),\\
\dar_H^G&:&D^b(kG)\lra D^b(kH),\\
\dar_K^H&:&D^b(kH)\lra D^b(kK)
\end{eqnarray*}
we get
$$\uar_H^G\circ\uar_K^H=\uar_K^G\textup{ and }\dar^H_K\circ\dar^G_H=\dar^G_K.$$
\end{Lemma}

\begin{proof}
This follows trivially from  the module case.
\end{proof}

Note that the notion of $\dar^G_H$-relative projectivity in $D^b(kG)$ corresponds to the similar
concept of relative projectivity with respect to a subalgebra as
developed in \cite[Section 2.1.1]{reptheobuch}. We shall need to extend the statements from there to
our more general situation.

\begin{Lemma}\label{verticesarepgroups}
Let $G$ be a finite group, and
let $k$ be a field of characteristic $p>0$. Let $D$ be a minimal subgroup of $G$ such that
the bounded complex of $kG$-modules
$M$ is $\dar^G_D$-relative projective. Then $D$ is a $p$-group.
\end{Lemma}

\begin{proof} Let $D\in Syl_p(G)$.
By Example~\ref{examplederivedcategoriesofgroups} every object $M$ of $D^b(kG)$ is
$\dar^G_D$-relative projective since $|G:D|$ is prime to $p$ by the definition of a Sylow subgroup.
If $M$ is $\dar^G_H$-relative projective, by Proposition~\ref{relativelyprojectivecharacterisation}
it is in $\add(\im\uar_H^G)$
and if $D'\in Syl_p(H)$, then $M$ is also
$\dar^H_{D'}$-relative projective, whence in $\add(\im\uar_{D'}^H)$ by
Proposition~\ref{relativelyprojectivecharacterisation} again. Therefore
$M$ is in $\add(\im\uar_{D'}^G)$, and therefore $\dar^G_{D'}$-relatively  projective,
again by Proposition~\ref{relativelyprojectivecharacterisation}.
\end{proof}

\begin{Def}
Let $G$ be a finite group, and
let $k$ be a field of characteristic $p>0$. Then, an indecomposable object
$M$ of $D^b(kG)$ has vertex $D$ if $M$ is relatively $kD$-projective, and
if $D$ is minimal with this property.
\end{Def}

\begin{Prop}\label{unicityofvertex}
The vertex $D$ of an indecomposable object $M$ of $D^b(kG)$ is a $p$-subgroup
of $G$, and $D$ is unique up to conjugacy.
\end{Prop}

\begin{proof} Using Lemma~\ref{verticesarepgroups}
we only need to show unicity up to conjugation.

The unicity part up to conjugation can be shown completely
analogous to the classical case.
Suppose that $M$ is a direct summand of $L\uar_K^G$ and of
$N\uar_H^G$ for two subgroups $H$ and $K$ of $G$ and two
indecomposable objects $L$ in $D^b(kK)$ and $N$ in $D^b(kH)$.
By Proposition~\ref{HigmansLemmaForAdjointFunctors}
we may suppose $L=M\dar^G_K$ and $N=M\dar^G_H$.
Then $M$ is a direct factor of
\begin{eqnarray*}
M\dar^G_H\uar_H^G\dar^G_K\uar_K^G&=&
\bigoplus_{KgH\in K\backslash G/H}{}^gM\dar^G_H\dar^H_{{}^gH\cap K}\uar_{{}^gH\cap K}^G\\
&=&\bigoplus_{KgH\in K\backslash G/H}{}^gM\dar^G_{{}^gH\cap K}\uar_{{}^gH\cap K}^G\\
\end{eqnarray*}
Using the Krull-Schmidt property, $M$ is a direct factor of
${}^gM\dar^G_{{}^gH\cap K}\uar_{{}^gH\cap K}^G$ for some $g$, and
since $K$ is minimal, there is $g\in G$ such that $^gH=K$.
\end{proof}

\begin{Rem}
The statements of the above results should remain true when we
replace this quite specific setting by
a Mackey functor with values in the functor category
between triangulated categories.
\end{Rem}

\begin{Lemma}\label{homologyisHprojective}
Let $G$ be a finite group, and let $k$ be a field of characteristic
$p>0$. Let $M$ be an indecomposable object of $D^b(kG)$. If $M$
is $kH$-projective, then each indecomposable direct factor of $H^n(M)$ for
all $n\in\N$ is relatively $H$-projective.
\end{Lemma}

\begin{proof}
By Proposition~\ref{relativelyprojectivecharacterisation} we see that
$M$ is relatively $D^b(kH)$-projective if and only if $M$ is a direct factor of
$L\uar_H^G$ for some $L$ in $D^b(kH)$. Since $\uar_H^G$ is exact, also
$H^n(L\uar_H^G)$ has a direct factor $H^n(M)$. However, $H^n(L\uar_H^G)=H^n(L)\uar_H^G$
and hence $H^n(M)$ is a direct factor of $H^n(L)\uar_H^G$. Hence, by
Higman's lemma \cite[Proposition 2.1.15]{reptheobuch} from modular representation theory,
each direct factor of $H^n(M)$  is relatively $H$-projective.
\end{proof}

\section{Localising on triangulated subcategories}

\label{VerdierlocalisingSection}

As in \cite{AuslanderKleiner} we consider the situation of three triangulated categories
with functors
$$
\xymatrix{
{\mathcal D}\ar@/^/[r]^{S'}&{\mathcal H}\ar@/^/[l]^{T'}\ar@/^/[r]^{S}&{\mathcal G}\ar@/^/[l]^T
}
$$
such that $(S,T)$ and $(S',T')$ are adjoint pairs. Let $\epsilon:\id_{\mathcal H}\lra TS$
be the unit of the adjunction $(S,T)$ and suppose that the unit is a split monomorphism. Hence
$$\xymatrix{\id_{\mathcal H}\ar[r]^\epsilon&TS\ar[r]^\pi&U\ar[r]^0&\id_{\mathcal H}}[1]$$
is a distinguished triangle of functors. In particular, $TS=\id_{\mathcal H}\oplus U$.

\begin{Rem}
Let $\mathcal T$ and $\mathcal S$ be triangulated categories, and let
$F:{\mathcal S}\lra{\mathcal T}$ be a triangle functor. Then with the
convention of Notion~\ref{inverseimageoffunctors} for any full
triangulated subcategory $\mathcal U$ of $\mathcal T$ the
category $F^{-1}({\mathcal U})$ is not necessarily triangulated.
However, if $\mathcal U$ is in addition closed under direct summands,
then this is true, as is shown in
Proposition~\ref{inverseimageofthick} below.
\end{Rem}

\begin{Prop}\label{inverseimageofthick}
Let $\mathcal T$ and $\mathcal S$ be triangulated categories, and let
$F:{\mathcal S}\lra{\mathcal T}$ be a triangle functor. Then for any thick
subcategory $\mathcal U$ of $\mathcal T$ the
category $F^{-1}({\mathcal U})$ is a triangulated subcategory of $\mathcal S$.
\end{Prop}

\begin{proof}
Let $X$ be an object of
${\mathcal S}$ such that $F(X)$ is a direct factor of the
object $U$ of $\mathcal U$. Hence, $F(X)\oplus U'=U$ for some object
$U'$  of $\mathcal T$. Since $\mathcal U$ is closed under direct factors,
$F(X)$ and $U'$ are actually already objects of $\mathcal U$.
Let  $X_1$ and $X_2$ be two objects of
${\mathcal S}$ such that $F(X_1)$ and $F(X_2)$ are objects of $\mathcal U$.
If
$$\xymatrix{X_1\ar[r]^\alpha&X_2\ar[r]&C(\alpha)\ar[r]&X_1[1]}$$
is a distinguished triangle in $\mathcal S$, since $F$ is a triangle functor, also
$$\xymatrix{FX_1\ar[r]^{F\alpha}&FX_2\ar[r]&FC(\alpha)\ar[r]&FX_1[1]}$$
is a distinguished triangle in $\mathcal T$, and hence $F(C(\alpha))\simeq C(F(\alpha))$.
Since $\mathcal U$ is triangulated $C(F(\alpha))$ is an object in $\mathcal U$,
and since $\mathcal U$ is closed under isomorphisms, $F(C(\alpha))$ is an
object of $\mathcal U$. Hence $C(\alpha)$ is an object of $F^{-1}({\mathcal U})$.
Therefore, $F^{-1}({\mathcal U})$ is a triangulated subcategory of $\mathcal S$.
\end{proof}

\begin{Lemma} \label{Zistriangulated}
Let $\mathcal Y$ be a full triangulated subcategory  of ${\mathcal H}$.
Then ${\mathcal Z}:=(US')^{-1}({\mathcal Y})$ satisfies
$S'({\mathcal Z})=S'({\mathcal D})\cap U^{-1}({\mathcal Y})$.
Moreover, $\mathcal Z$ is triangulated if $\mathcal Y$ is thick.
\end{Lemma}

\begin{proof} By definition
$S'({\mathcal Z})$ is the full subcategory of $\mathcal H$
formed by objects $S'M$ such that $US'M\in \add({\mathcal Y})$. Hence
$S'({\mathcal Z})$ is contained in $S'({\mathcal D})\cap U^{-1}({\mathcal Y})$.
Moreover, an object $X$ in $S'({\mathcal D})\cap U^{-1}({\mathcal Y})$ is
an object of the form $S'M$, since  $X\in S'({\mathcal D})$ and such that $US'M\in \add({\mathcal Y})$
since $X\in U^{-1}({\mathcal Y})$. The rest follows from Proposition~\ref{inverseimageofthick}.
\end{proof}

\begin{Rem}\label{Remarkdistinguishthenotations}
We remind the reader that we have two different localisation or quotient
constructions of a triangulated category $\mathcal T$ by a triangulated
subcategory $\mathcal U$ (cf Notion~\ref{inverseimageoffunctors}).
\begin{itemize}
\item
First we have the additive quotient, denoted traditionally ${\mathcal S}/{\mathcal U}$
having the same objects as $\mathcal S$ but we consider morphisms between
two objects as residue classes of morphisms in $\mathcal T$ modulo those factoring
through an object of $\mathcal U$.
\item
Second, the Verdier localisation \cite{SGA412,Verdier} which we
denote by ${\mathcal S}_{\mathcal U}$.
In the literature the Verdier localisation is often denoted by
${\mathcal S}/{\mathcal U}$. In order to distinguish from the additive
quotient we decided to use the symbol ${\mathcal S}_{\mathcal U}$, contrary to
the established convention in the literature.
\end{itemize}
\end{Rem}

\begin{Lemma}\label{naturalquotientinduces}
Let $\mathcal S$ be a triangulated category, and let $\mathcal U$ be a
thick subcategory of $\mathcal S$. Then there is a unique and natural functor
${\mathcal S}/{\mathcal U}\stackrel{L_{\mathcal U}}\lra {\mathcal S}_{\mathcal U}$ making the
diagram
$$
\xymatrix{&{\mathcal S}\ar[dl]_{Q_{\mathcal U}}\ar[dr]^{V_{\mathcal U}}\\
{\mathcal S}/{\mathcal U}\ar@{-->}[rr]_{L_{\mathcal U}}&&{\mathcal S}_{\mathcal U}
}
$$
commutative. Here we denote
${\mathcal S}\stackrel {Q_{\mathcal U}}{\lra} {\mathcal S}/{\mathcal U}$ and
${\mathcal S}\stackrel {V_{\mathcal U}}\lra {\mathcal S}_{\mathcal U}$ the canonical functors
given by the respective universal properties.
\end{Lemma}

\begin{proof}
The proof is implicit in \cite[Proposition 1.3]{Rickardstable}.
The statement follows from the well-known fact that
for any additive category
$\mathcal A$ and any additive functor $F:{\mathcal S}\lra{\mathcal A}$ such that
$F(U)=0$ for any
object $U$ of $\mathcal U$, there is a unique additive functor
$F^*:{\mathcal S}/{\mathcal U}\lra \mathcal A$ with $F=F^*\circ Q_{\mathcal U}$.
For $\mathcal A={\mathcal S}_{\mathcal U}$, we observe that
$F=V_{\mathcal U}:{\mathcal S}\rightarrow {\mathcal S}_{\mathcal U}$ is additive
with $F(U)=0$ for any object $U$ of $\mathcal U$. Indeed, a morphism in the
localisation becomes invertible if its cone is in $\mathcal U$. Hence, for an
object $U$ of $\mathcal U$ the cone of the zero morphism on $U$ is in $\mathcal U$.
Therefore the image $V_{\mathcal U}(0_U)$ of the zero morphism $0_U$ on $U$ in
the localisation is invertible in ${\mathcal S}_{\mathcal U}$.
The only object with invertible zero endomorphism is the zero object in
${\mathcal S}_{\mathcal U}$. This proves the statement.
%
\end{proof}

\begin{Rem}
We need to recall from Verdier~\cite[Chapitre II, Section 2.1, 2.2]{Verdier},
or alternatively from Stacks project \cite[Part 1, Chapter 13, Section 13.6]{Stacks}, some
properties of Verdier localisation.
If $F:{\mathcal S}\lra{\mathcal T}$ is a triangle functor between triangulated categories,
then the full subcategory $\ker(F)$ of $\mathcal S$ generated by those objects $X$ of $\mathcal S$
such that $F(X)=0$ is thick.
If $\mathcal U$ is a full triangulated subcategory of some triangulated category
$\mathcal T$, then the Verdier localisation defined by inverting all morphisms
$f$ in $\mathcal T$
with cone in $\mathcal U$ is triangulated, and there is a
canonical functor $V_{\mathcal U}:{\mathcal T}\lra{\mathcal T}_{\mathcal U}$
with ${\mathcal U}$ is a full triangulated subcategory of $\ker(V_{\mathcal U})$.
Moreover, $\ker(V_{\mathcal U})$ is thick,
namely the smallest thick subcategory of
$\mathcal T$ containing $\mathcal U$, the thickening $\thick({\mathcal U})$.
\end{Rem}

\begin{Rem}
We see that even if $\mathcal Y$ is a thick subcategory of the triangulated category
$\mathcal C$ and if $H$ is a triangle functor ${\mathcal C}\ra{\mathcal D}$
for some triangulated category $\mathcal D$, then $H({\mathcal Y})$ is triangulated, but
is not thick anymore in general.
The Verdier localisation of $\mathcal D$ at $H({\mathcal Y})$ has good
properties with respect to thick subcategories.
Since $\ker(V_{H({\mathcal Y})})=\thick(H({\mathcal Y}))$, we need to consider
$\thick (H({\mathcal Y}))$.
Then
$$\ker(\xymatrix{{\mathcal T}\ar[rr]^-{V_{\thick(H{\mathcal Y})}}&&
{\mathcal T}_{\thick (H({\mathcal Y}))}})=\thick (H({\mathcal Y}))=
\ker(\xymatrix{{\mathcal T}\ar[r]^-{V_{H\mathcal Y}}& {\mathcal T}_{H\mathcal Y}}).$$
\end{Rem}

\begin{Lemma}\label{extendtolocaliseatthick}
Let $\mathcal C$ and $\mathcal D$ be triangulated categories, let $\mathcal Y$ be a
 subcategory of $\mathcal C$, and let $H:{\mathcal C}\ra{\mathcal D}$ be a triangle functor.
Then $H$ extends to a unique functor
${\mathcal C}/{\mathcal Y}\ra{\mathcal D}_{\thick H({\mathcal Y})}$, also denoted by $H$,
such that $H\circ Q_{\mathcal Y}=V_{\thick H\mathcal Y}\circ H$, i.e. making the diagram
$$
\xymatrix{
{\mathcal C}\ar[rr]^H\ar[d]_{Q_{\mathcal Y}}&&
{\mathcal D}\ar[d]^{V_{\thick H\mathcal Y}}\ar[dl]_{Q_{H\mathcal Y}}\\
{\mathcal C}/{\mathcal Y}\ar[r]^H\ar@/_1pc/[rr]_H&
{\mathcal D}/H{\mathcal Y}\ar[r]&{\mathcal D}_{\thick H\mathcal Y}
}
$$
commutative. The functor ${\mathcal D}/{H\mathcal Y}\lra {\mathcal D}/{\thick H\mathcal Y}$
combined with the functor
$L_{\thick H\mathcal Y}:{\mathcal D}/{\thick H\mathcal Y}\lra
{\mathcal D}_{\thick H\mathcal Y}$ make the right triangle of the above diagram commutative.
\end{Lemma}

\begin{proof}
Consider
$$\xymatrix{{\mathcal C}\ar[r]^{H}&{\mathcal D}\ar[r]^-{V_{\thick H\mathcal Y}}&
{\mathcal D}_{\thick H{\mathcal Y}}.}$$
Then for all objects $X$ of $\mathcal Y$ we get $(V_{\thick H\mathcal Y}\circ H)(X)=0$.
Likewise consider
$$\xymatrix{{\mathcal C}\ar[r]^{H}&{\mathcal D}\ar[r]^-{Q_{H\mathcal Y}}&
{\mathcal D}/{H{\mathcal Y}}.}$$
Again, for all objects $X$ of $\mathcal Y$ we get $(Q_{H\mathcal Y}\circ H)(X)=0$.
Hence there is a unique functor
${\mathcal C}/{\mathcal Y}\stackrel{H_1}{\lra}{\mathcal D}_{\thick H({\mathcal Y})}$
satisfying $H_1\circ Q_{\mathcal Y}=V_{\thick H\mathcal Y}\circ H$
and a unique functor ${\mathcal C}/{\mathcal Y}\stackrel{H_2}{\lra}{\mathcal D}/{H({\mathcal Y})}$
satisfying $H_2\circ Q_{\mathcal Y}=Q_{H\mathcal Y}\circ H$.
Moreover, by the universal property of ${\mathcal D}/{H\mathcal Y}$
the functor $H_1$ factors through $H_2$ and through $L_{\thick H\mathcal Y}$.
\end{proof}

\begin{Lemma}\label{refereeslemma}
Let $F:{\mathcal C}\lra {\mathcal D}$ be a triangle functor between triangulated categories,
let $\mathcal X$ be a full subcategory of $\mathcal C$, and let $\mathcal Y$ be a full
subcategory of $\mathcal D$. If $F({\mathcal X})\subseteq \mathcal Y$, then there exists
a unique additive functor ${\mathcal C}_{\thick\mathcal X}\lra {\mathcal D}_{\thick\mathcal Y}$,
still denoted by $F$, making the following diagram commutative.
$$\xymatrix{
{\mathcal C}/{\mathcal X}\ar[d]_F\ar[r]&
{\mathcal C}/{\thick\mathcal X}\ar[d]_F\ar[r]^{L_{\thick\mathcal X}}&
{\mathcal C}_{\thick\mathcal X}\ar[d]_F\\
{\mathcal D}/{\mathcal Y}\ar[r]
&{\mathcal D}/{\thick\mathcal Y} \ar[r]_{L_{\thick\mathcal Y}}&{\mathcal D}_{\thick\mathcal Y}
}$$
\end{Lemma}

\begin{proof}
Since $F$ is a triangle functor,
$$F(\thick{\mathcal X})\subseteq \thick F{\mathcal X}\subseteq \thick{\mathcal Y}$$
and the left square is commutative.
By Lemma~\ref{extendtolocaliseatthick} the right square is commutative as well.
\end{proof}

For a triangulated category $\mathcal T$ and a full subcategory $\mathcal X$
there is a natural functor ${\mathcal T}/{\mathcal X}\lra {\mathcal T}/\thick{\mathcal X}$.
For simplicity the composition
$$\xymatrix{{\mathcal T}/{\mathcal X}\ar[r]&
{\mathcal T}/\thick{\mathcal X}\ar[r]^-{L_{\thick\mathcal X}}&{\mathcal T}_{\thick\mathcal X}
}$$
is also denote by $L_{\thick\mathcal X}$.

\begin{Prop} \label{functorsfrommodtoVerdier}
Let ${\mathcal D}$, ${\mathcal H}$, ${\mathcal G}$ be three triangulated categories
and triangle functors $S,S',T,T'$
$$
\xymatrix{
{\mathcal D}\ar@/^/[r]^{S'}&{\mathcal H}\ar@/^/[l]^{T'}\ar@/^/[r]^{S}&{\mathcal G}\ar@/^/[l]^T
}
$$
so that $(S,T)$ and $(S',T')$ are adjoint pairs. Let $\epsilon:\id_{\mathcal H}\lra TS$
be the unit of the adjunction $(S,T)$. Assume that there is an
endofunctor $U$ of $\mathcal H$ such that $TS=\id_{\mathcal H}\oplus U$, denote by
$p_1:TS\lra \id_{\mathcal H}$
the projection, and suppose that $p_1\circ\epsilon$ is an isomorphism.

Let ${\mathcal Y}$ be a thick subcategory of $\mathcal H$, and suppose that
each object of $TSS'T'{\mathcal Y}$ is a direct factor of an object of
${\mathcal Y}$. Then
$S$ and $T$ induce triangle functors
$${\mathcal H}_{\thick S'T'{\mathcal Y}}\stackrel{S}{\lra}
{\mathcal G}_{\thick SS'T'{\mathcal Y}}\mbox{ and }
{\mathcal G}_{\thick SS'T'{\mathcal Y}}\stackrel{T}{\lra}{\mathcal H}_{\mathcal Y}$$
making the diagram
$$
\xymatrix{
{\mathcal H}/{\thick S'T'{\mathcal Y}}\ar[r]^-{S}\ar[d]_{L_{\thick S'T'{\mathcal Y}}}&
{\mathcal G}/{\thick SS'T'{\mathcal Y}}\ar[d]^{L_{\thick SS'T'{\mathcal Y}}}\ar[r]^-T&
{\mathcal H}/{{\mathcal Y}}\ar[d]^{L_{\mathcal Y}}\\
{\mathcal H}_{\thick S'T'{\mathcal Y}}\ar[r]^-{S}&
{\mathcal G}_{\thick SS'T'{\mathcal Y}}\ar[r]^-T&{\mathcal H}_{{\mathcal Y}}
}
$$
commutative.
\end{Prop}

\begin{proof}
The existence of the functors in the left square and the commutativity of the left square
follow from Lemma~\ref{extendtolocaliseatthick}. From Lemma~\ref{extendtolocaliseatthick}
we get natural functors giving a commutative diagram
$$
\xymatrix{
{\mathcal H}/{S'T'{\mathcal Y}}\ar[r]^-{S}\ar[d]_{L_{\thick S'T'{\mathcal Y}}}&
{\mathcal G}/{SS'T'{\mathcal Y}}\ar[d]^{L_{\thick SS'T'{\mathcal Y}}}\ar[r]^T&
{\mathcal H}/{TSS'T'{\mathcal Y}}\ar[d]^{L_{\thick SS'T'\mathcal Y}}\\
{\mathcal H}_{\thick S'T'{\mathcal Y}}\ar[r]^-{S}&
{\mathcal G}_{\thick SS'T'{\mathcal Y}}\ar[r]^T&{\mathcal H}_{\thick TSS'T'{\mathcal Y}}
}
$$
For ${\mathcal X}=SS'T'\mathcal Y$,  Theorem~\ref{Greenforadjoints} shows that
$T({\mathcal X})\subseteq \mathcal Y$.
Using Lemma~\ref{refereeslemma} and the fact that $\mathcal Y$ is thick,
and therefore $\thick\mathcal Y=\mathcal Y$, we obtain a commutative diagram
$$
\xymatrix{
{\mathcal H}/{S'T'{\mathcal Y}}\ar[r]^-{S}\ar[d]_{L_{\thick S'T'{\mathcal Y}}}&
{\mathcal G}/{SS'T'{\mathcal Y}}\ar[d]^{L_{\thick SS'T'{\mathcal Y}}}\ar[r]^-T&
{\mathcal H}/{\mathcal Y}\ar[d]^{L_{{\mathcal Y}}}\\
{\mathcal H}_{\thick S'T'{\mathcal Y}}\ar[r]^-{S}&
{\mathcal G}_{\thick SS'T'{\mathcal Y}}\ar[r]^-T&{\mathcal H}_{{\mathcal Y}}
}
$$
as requested.

The fact that
$${\mathcal H}_{\thick S'T'{\mathcal Y}}\stackrel{S}{\lra} {\mathcal G}_{\thick SS'T'{\mathcal Y}}
\mbox{ and }
{\mathcal G}_{\thick SS'T'{\mathcal Y}}\stackrel{T}{\lra}{\mathcal H}_{\mathcal Y}$$
are triangle functors comes from the universal property of the Verdier localisation
(cf \cite[Chapitre II, Theor\`eme 2.2.6]{Verdier}).
%
\end{proof}

\begin{Cor}
Let ${\mathcal D}$, ${\mathcal H}$, ${\mathcal G}$ be three triangulated categories
and triangle functors $S,S',T,T'$
$$
\xymatrix{
{\mathcal D}\ar@/^/[r]^{S'}&{\mathcal H}\ar@/^/[l]^{T'}\ar@/^/[r]^{S}&{\mathcal G}\ar@/^/[l]^T
}
$$
so that $(S,T)$ and $(S',T')$ are adjoint pairs. Let $\epsilon:\id_{\mathcal H}\lra TS$
be the unit of the adjunction $(S,T)$. Assume that there is an
endofunctor $U$ of $\mathcal H$ such that $TS=\id_{\mathcal H}\oplus U$, denote by
$p_1:TS\lra \id_{\mathcal H}$
the projection, and suppose that $p_1\circ\epsilon$ is an isomorphism.

Let ${\mathcal Y}$ be a thick subcategory of $\mathcal H$, and suppose that
each object of $TSS'T'{\mathcal Y}$ is a direct factor of an object of
${\mathcal Y}$.
Then we have a commutative diagram
$$
\xymatrix{
{\mathcal H}_{\thick S'T'{\mathcal Y}}\ar[r]|-S\ar[d]_{\text{can}}&
{\mathcal G}_{\thick SS'T'{\mathcal Y}}\ar[dl]|-T\\
{\mathcal H}_{{\mathcal Y}}\\
&&{\mathcal H}/{ S'T'{\mathcal Y}}\ar[r]|S\ar[d]_{\text{can}}\ar@/_/[uull]|{L_{\thick S'T'{\mathcal Y}}}&
{\mathcal G}/{ SS'T'{\mathcal Y}}\ar[dl]|T\ar@/_/[uull]|{L_{\thick SS'T'{\mathcal Y}}}\\
&&{\mathcal H}/{{\mathcal Y}}\ar@/_/[uull]|{L_{\mathcal Y}}
}
$$
\end{Cor}

\begin{proof}
Indeed, since each object of $TSS'T'{\mathcal Y}$ is a direct factor of an object of
${\mathcal Y}$, each object of $S'T'\mathcal Y$ is a direct factor of an object of
$\mathcal Y$.
Hence, there is a natural functor $\text{can}$ as indicated. The rest of the statement
is an immediate consequence of Proposition~\ref{functorsfrommodtoVerdier}.
\end{proof}

\begin{Cor}\label{functorsextendtoVerdierlocalisations}
Let ${\mathcal D}$, ${\mathcal H}$, ${\mathcal G}$ be three triangulated categories
and triangle functors $S,S',T,T'$
$$
\xymatrix{
{\mathcal D}\ar@/^/[r]^{S'}&{\mathcal H}\ar@/^/[l]^{T'}\ar@/^/[r]^{S}&{\mathcal G}\ar@/^/[l]^T
}
$$
so that $(S,T)$ and $(S',T')$ are adjoint pairs. Let $\epsilon:\id_{\mathcal H}\lra TS$
be the unit of the adjunction $(S,T)$. Assume that there is an
endofunctor $U$ of $\mathcal H$ such that $TS=\id_{\mathcal H}\oplus U$, denote by
$p_1:TS\lra \id_{\mathcal H}$
the projection, and suppose that $p_1\circ\epsilon$ is an isomorphism.

Let ${\mathcal Y}$ be a thick subcategory of $\mathcal H$,
put ${\mathcal Z}:= (US')^{-1}({\mathcal Y})$, and suppose that
each object of $TSS'T'{\mathcal Y}$ is a direct factor of an object of
${\mathcal Y}$.
Then the restriction of $S$ to the subcategory ${\add S'{\mathcal Z}}/{ S'T'{\mathcal Y}}$ and
the restriction of $T$ to
the subcategory ${\add SS'{\mathcal Z}}/{ SS'T'{\mathcal Y}}$ are equivalences and gives
a commutative diagram
$$
\xymatrix{
{\thick S'{\mathcal Z}}_{\thick S'T'{\mathcal Y}}\ar[r]|-S\ar[d]_{\text{can}}&
{\thick SS'{\mathcal Z}}_{\thick SS'T'{\mathcal Y}}\ar[dl]|-T\\
{\thick S'{\mathcal Z}}_{{\mathcal Y}}\\
&&{\add S'{\mathcal Z}}/{ S'T'{\mathcal Y}}\ar[r]|-S\ar[d]_{\text{can}}
\ar@/_/[uull]|{L_{\thick S'T'{\mathcal Y}}}&
{\add SS'{\mathcal Z}}/{ SS'T'{\mathcal Y}}\ar[dl]|-T\ar@/_/[uull]|{L_{\thick SS'T'{\mathcal Y}}}\\
&&{\add TSS'{\mathcal Z}}/{{\mathcal Y}}\ar@/_/[uull]|{L_{\mathcal Y}}
}
$$
where the lower triangle consists of equivalences.
\end{Cor}

\begin{proof}
The fact that the lower triangle exists and is commutative follows from
Theorem~\ref{Greenforadjoints}. Since for any subcategory $\mathcal X$ of $\mathcal T$
we get that $\add{\mathcal X}$ is a full subcategory of $\thick \mathcal X$,
we have a commutative diagram
$$
\xymatrix{
({\thick S'{\mathcal Z}})/{S'T'{\mathcal Y}}\ar[r]|-S\ar[d]_{\text{can}}&({\thick SS'{\mathcal Z}})/{ SS'T'{\mathcal Y}}\ar[dl]|-T\\
({\thick S'{\mathcal Z}})/{{\mathcal Y}}\\
&&{\add S'{\mathcal Z}}/{ S'T'{\mathcal Y}}\ar[r]|-S\ar[d]_{\text{can}}
\ar@/_/[uull]&{\add SS'{\mathcal Z}}/{ SS'T'{\mathcal Y}}\ar[dl]|-T\ar@/_/[uull]\\
&&{\add TSS'{\mathcal Z}}/{{\mathcal Y}}\ar@/_/[uull]
}
$$
By Lemma~\ref{extendtolocaliseatthick} we obtain a commutative diagram
$$
\xymatrix{
{\thick S'{\mathcal Z}}_{\thick S'T'{\mathcal Y}}\ar[r]|-S\ar[d]_{\text{can}}&
{\thick SS'{\mathcal Z}}_{\thick SS'T'{\mathcal Y}}\ar[dl]|-T\\
{\thick S'{\mathcal Z}}_{{\mathcal Y}}\\
&&({\thick S'{\mathcal Z}})/{S'T'{\mathcal Y}}\ar[r]|-S\ar[d]_{\text{can}}
\ar@/_/[uull]|{L_{\thick S'T'{\mathcal Y}}}&
({\thick SS'{\mathcal Z}})/{ SS'T'{\mathcal Y}}\ar[dl]|-T\ar@/_/[uull]|{L_{\thick SS'T'{\mathcal Y}}}\\
&&({\thick S'{\mathcal Z}})/{{\mathcal Y}}\ar@/_/[uull]|{L_{\mathcal Y}}
}
$$
Composition of the two diagrams yields the statement.
\end{proof}

\begin{Prop}\label{thickeningnumeratorinordinaryquotient}
Let ${\mathcal T}$ and ${\mathcal U}$ be two triangulated categories, let
$F:{\mathcal T}\lra {\mathcal U}$ be a triangle functor,
let ${\mathcal X}$ be a full
additive subcategory of $\mathcal T$ and let ${\mathcal Y}$ be a full additive subcategory
of $\mathcal X$.
Then the restriction of $F$ to
${\mathcal X}/{\mathcal Y}\stackrel{F_{\mathcal X}}\lra (\add F{\mathcal X})/F{\mathcal Y}$
extends to a functor
$(\thick {\mathcal X})/{\mathcal Y}\stackrel{F_{\thick\mathcal X}}\lra (\thick F{\mathcal X})/F{\mathcal Y}$ such that $F_{\thick\mathcal X}$ coincides with $F_{\mathcal X}$ on the subcategory ${\mathcal X}/{\mathcal Y}$.
\end{Prop}

\begin{proof}
Let $\mathcal A$ and $\mathcal B$ be full subcategories of a triangulated
category $\mathcal V$, then as in \cite{BBD} we denote by ${\mathcal A}\ast{\mathcal B}$
the full subcategory of $\mathcal V$ generated by $C(t)[-1]$ where
$$\xymatrix{
A\ar[r]^-f&C(t)[-1]\ar[r]^-s&B\ar[r]^-t&A[1]
}$$
is a distinguished triangle, and where $A$ is an object of $\mathcal A$, $B$ is an object of
$\mathcal B$, and $t\in{\mathcal T}(B,A[1])$.
Since  $F:{\mathcal T}\lra{\mathcal U}$ is a triangle functor, $F$ sends distinguished
triangles to distinguished triangle. Therefore $F({\mathcal A}\ast{\mathcal B})$ is a
subcategory of
$F({\mathcal A})\ast F({\mathcal B})$.
Hence $F$ induces a functor
$$\xymatrix{
{\mathcal A}\ast{\mathcal B}\ar[r]^-F&(F{\mathcal A})\ast(F{\mathcal B}).
}$$

Let $X_1\in\add F{\mathcal A}$ and $X_2\in\add F{\mathcal B}$.
Then for any $t\in{\mathcal U}(X_2,X_1[1])$ we get
$$C(t)[-1]\in\add(F({\mathcal A})\ast F({\mathcal B})).$$
Indeed, denote by
$$\xymatrix{X_1\ar[r]^-f&C(t)[-1]\ar[r]^-s&X_2\ar[r]^-t&X_1[1]}$$
the distinguished triangle given by $t$.
Let $X_1'$ and $X_2'$ be objects of $\mathcal U$ such that
$X_1\oplus X_1'\in F({\mathcal A})$ and $X_2\oplus X_2'\in F({\mathcal B})$.
Then
$$
\xymatrix{
X_1\oplus X_1'\ar[rr]^-{\left(\begin{array}{cc}f&0\\0&\id_{X_1'}\\0&0\end{array}\right)}&&
C(t)[-1]\oplus X_1'\oplus X_2'\ar[rr]^-{\left(\begin{array}{ccc}s&0&0\\0&0&\id_{X_2'}\end{array}\right)}&&
X_2\oplus X_2'\ar[rr]^-{\left(\begin{array}{cc}t&0\\0&0\end{array}\right)}&&(X_1\oplus X_1')[1]}
$$
is a distinguished triangle. Hence $C(t)[-1]\in\add(F({\mathcal A})\ast F({\mathcal B}))$.
This shows
$$\add(F({\mathcal A}))\ast \add(F({\mathcal B}))\subseteq\add(F({\mathcal A})\ast F({\mathcal B})),$$
and therefore
$$\add\left(\add(F({\mathcal A}))\ast \add(F({\mathcal B}))\right)=\add(F({\mathcal A})\ast F({\mathcal B})).$$
If we define $({\mathcal Z})_n:=({\mathcal Z})_{n-1}\ast {\mathcal Z}$ for any
subcategory $\mathcal Z$ of $\mathcal U$, and
$({\mathcal Z})_1:={\mathcal Z}$, then
$$\thick(F({\mathcal X}))=\bigcup_{n\in\N}\add((\add(F({\mathcal X})))_n)=
\bigcup_{n\in\N}\add((F({\mathcal X}))_n)=\add(\bigcup_{n\in\N}(F({\mathcal X}))_n).$$
Now,
$$\thick{\mathcal X}=\bigcup_{n\in\N}\add{\mathcal X}_n=\add\left(\bigcup_{n\in\N}{\mathcal X}_n\right)$$
and $F({\mathcal X}_n)\subseteq (F{\mathcal X})_n$.
Hence
\begin{eqnarray*}
F(\thick{\mathcal X})/F{\mathcal Y}&=&F(\add\left(\bigcup_{n\in\N}{\mathcal X}_n\right))/F{\mathcal Y}\\
&\subseteq&\add\left(\bigcup_{n\in\N}F({\mathcal X}_n)\right)/F{\mathcal Y}\\
&\subseteq&\add\left(\bigcup_{n\in\N}(F({\mathcal X})_n)\right)/F{\mathcal Y}\\
&=&\thick(F({\mathcal X}))/F{\mathcal Y}
\end{eqnarray*}
Therefore $F$ extends to a functor
$$\thick({\mathcal X})/{\mathcal Y}\stackrel{F_{\thick\mathcal X}}\lra \thick(F({\mathcal X}))/F{\mathcal Y}.$$
By construction the restriction of $F_{\thick{\mathcal X}}$ to ${\mathcal X}/{\mathcal Y}$
coincides with $F_{\mathcal X}$.
\end{proof}

\begin{Rem}
If in Proposition~\ref{thickeningnumeratorinordinaryquotient}
the functor $F$ induces an equivalence
${\mathcal X}/{\mathcal Y}\lra (\add F{\mathcal X})/F{\mathcal Y}$, then there is no reason
why this should imply an equivalence
$$\thick({\mathcal X})/{\mathcal Y}\lra \thick(F({\mathcal X}))/F{\mathcal Y}.$$
\end{Rem}

\begin{Lemma}\label{ordinaryquotientindependentofthickening}
Let $\mathcal Y$ be a subcategory of $\mathcal X$ admitting finite direct sums. Then the natural projection
${\mathcal X}/{\mathcal Y} \lra {\mathcal X}/(\add{\mathcal Y})$ is an equivalence of categories.
\end{Lemma}

\begin{proof}
Since ${\mathcal Y}$ is a subcategory of $\add{\mathcal Y}$, if a morphism
$f$ factors through an object of ${\mathcal Y}$, it factors also through an object of
$\add{\mathcal Y}$. Hence, the natural projection is well-defined and full.
If $f$ factors through an object $X$ of $\add{\mathcal Y}$, then there is an object
$X'$ of $\add{\mathcal Y}$, such that $X\oplus X'$ is an object of ${\mathcal Y}$.
Extending by the zero morphism to and from $X'$, hence $f$ factors also through the
object $X\oplus X'$ of $\mathcal Y$. This shows that the natural projection is
faithful as well.

From the above it also follows that the natural projection is dense,
since the objects of both quotient categories coincide,
and the natural projection is the identity on objects.
\end{proof}

\begin{Prop}\label{ordinaryquotientisogivesVerdierquotientiso}
Let ${\mathcal T}$ and ${\mathcal U}$ be triangulated categories, let ${\mathcal Y}$
be a subcategory of ${\mathcal T}$, and let $F:{\mathcal T}\lra{\mathcal U}$ be a triangle functor.
Suppose that $F$ induces an equivalence
$$F_Q:{\mathcal T}/(\thick{\mathcal Y})\lra {\mathcal U}/(\thick F({\mathcal Y})).$$
Then $F$ induces a dense and full triangle functor
$$F_V:{\mathcal T}_{(\thick{\mathcal Y})}\lra {\mathcal U}_{(\thick F({\mathcal Y}))}.$$
If in addition $F(\thick{\mathcal Y})$ is thick in $\mathcal U$
then $F_V$ is a triangle equivalence.
\end{Prop}

\begin{proof} The functor $F_V$ exists by the universal property of the Verdier localisation
\cite[Chapitre II, Corollaire 2.2.11.c]{Verdier}.

{\bf We shall now show that $F_V$ is dense.}
The objects of ${\mathcal U}/(\thick F({\mathcal Y}))$ coincides with the
objects of ${\mathcal U}_{(\thick F({\mathcal Y}))}$, since they both coincides
with the objects of $\mathcal U$. By hypothesis, for every object $U$ of $\mathcal U$
there is an object $T$ of $\mathcal T$, and $f\in{\mathcal U}(FT,U)$ as well as
$g\in{\mathcal U}(U,FT)$ such that $g\circ f-\id_{FT}$ factors through
an object $Y'$ of $\thick(F{\mathcal Y})$ and $f\circ g-\id_U$ factors through
an object $Y$ of $\thick({F\mathcal Y})$. Hence, applying $L_{\thick F{\mathcal Y}}$
to these equations, and
observing that $L_{\thick{F\mathcal Y}}(Y)=0$, respectively $L_{\thick F{\mathcal Y}}(Y')=0$
for all objects $Y$ in $\thick(F{\mathcal Y})$, respectively  all objects $Y'$ in $\thick (F{\mathcal Y})$,
we get that the image of $f$ in the Verdier localisation is an isomorphism.
Hence $F_V$ is dense.

\medskip

{\bf We will show now that $F_V$ is full.}

{\em First step:}
Let $f\in{\mathcal U}(FZ,FX)$. Since
$F_Q$ is full, there is $f'\in {\mathcal T}(Z,X)$ such that $f-Ff'$ factors through
an object $M$ of $\thick(F{\mathcal Y})$.
Hence there is $g\in{\mathcal U}(M,X)$ and $h\in{\mathcal U}(Z,M)$
with $f-Ff'=g\circ h$ in $\mathcal U$.
We denote by $(1,f)$ the morphism
represented by the diagram
$\xymatrix{FZ&\ar[l]_-{\id_{FZ}}FZ\ar[r]^f&FX}$. Then
$(1,f-Ff')=(1,g)\circ (1,h)$ in ${\mathcal U}_{(\thick F({\mathcal Y}))}$.
Since $M\simeq 0$ in ${\mathcal U}_{(\thick F({\mathcal Y}))}$, we get
$$(1,f)-(1,Ff')=(1,f-Ff')=0$$ in ${\mathcal U}_{(\thick F({\mathcal Y}))}$
and therefore $$(1,f)=(1,Ff')=F_V(1,f').$$
Hence $f=F_Qf'$ in ${\mathcal U}/\thick F{\mathcal Y}$ implies
$(1,f)=F_V(1,f')$ in ${\mathcal U}_{(\thick F({\mathcal Y}))}$.

{\em Second step:}
Let $\xymatrix{FX&\hat Z\ar[l]_-s\ar[r]^-f&FY}$ represent a morphism
in ${\mathcal U}_{(\thick F({\mathcal Y}))}(F_VX, F_VY)$. Since $F_V$ is
dense, we may suppose that $\hat Z=FZ$ for some object $Z$ of $\mathcal U$.
Then $s=F_Qs'$ for some $s'\in{\mathcal T}(Z,X)$, and $f=F_Qf'$ for some
$f'\in{\mathcal T}(Z,Y)$, giving $(1,s)=F_V(1,s')$, and $(1,f)=F_V(1,f')$ by the first step.
Now, $$F(\cone(s'))\simeq \cone(F(s'))\simeq \cone(s)\in\thick(F({\mathcal Y}))$$
and therefore $F(\cone(s'))=0$ in ${\mathcal U}/(\thick F{\mathcal Y})$. Since $F$ induces
an equivalence
$${\mathcal U}/(\thick F{\mathcal Y})\simeq {\mathcal T}/(\thick {\mathcal Y}),$$
we get $\cone(s')=0$ in ${\mathcal T}/(\thick {\mathcal Y})$, which shows
that $\cone(s')\in\thick{\mathcal Y}$.
Hence $\xymatrix{X&Z\ar[l]_-{s'}\ar[r]^-{f'}&Y}$
maps to $\xymatrix{FX&FZ\ar[l]_-s\ar[r]^-f&FY}$.
Therefore $F_V$ is full.

\medskip

We now assume that in addition $F(\thick{\mathcal Y})$ is thick in $\mathcal U$.

{\bf We need to show that $F_V$ is faithful.}
Since $F(\thick{\mathcal Y})$ is thick
in $\mathcal U$, we get  $\thick(F({\mathcal Y}))=F(\thick{\mathcal Y})$.
Hence, using the notation from \cite{Objective}, and using that $F_Q$ is an equivalence,
$\ker(F_V)=0$.
By a result of Rickard~\cite[first paragraph on page 446]{RickardMoritatheorem}
or Ringel and Zhang~\cite[Proposition 3.1, Proposition 3.3, Theorem 1.1]{Objective}
we see that $F_V$ is faithful.
This finishes the proof.
\end{proof}

\begin{Theorem} (Green correspondence for triangulated categories)
\label{Greencorrespondencefortriangulated}
Let ${\mathcal D}$, ${\mathcal H}$, ${\mathcal G}$ be three triangulated categories
and let $S,S',T,T'$ be triangle functors
$$
\xymatrix{
{\mathcal D}\ar@/^/[r]^{S'}&{\mathcal H}\ar@/^/[l]^{T'}\ar@/^/[r]^{S}&{\mathcal G}\ar@/^/[l]^T
}
$$
such that $(S,T)$ and $(S',T')$ are adjoint pairs. Let $\epsilon:\id_{\mathcal H}\lra TS$
be the unit of the adjunction $(S,T)$. Assume that there is an
endofunctor $U$ of $\mathcal H$ such that $TS=\id_{\mathcal H}\oplus U$, denote by
$p_1:TS\lra \id_{\mathcal H}$
the projection, and suppose that $p_1\circ\epsilon$ is an isomorphism.

Let ${\mathcal Y}$ be a thick subcategory of $\mathcal H$,
put ${\mathcal Z}:= (US')^{-1}({\mathcal Y})$, and suppose that
each object of $TSS'T'{\mathcal Y}$ is a direct factor of an object of
${\mathcal Y}$.
\begin{enumerate}
\item\label{(1)}
Then $S$ and $T$ induce triangle functors $S_Z$ and $T_Z$ fitting into the commutative diagram
$$
\xymatrix{
({\thick (S'{\mathcal Z})})_{(\thick (S'T'{\mathcal Y}))}\ar[r]|-{S_Z}\ar[d]_{\text{can}}&
({\thick (SS'{\mathcal Z})})_{(\thick (SS'T'{\mathcal Y}))}\ar[dl]|-{T_Z}\\
(\thick (S'{\mathcal Z}))_{{\mathcal Y}}}
$$
of Verdier localisations.
\item\label{(2)}
There is an additive functor $S_\thick$, induced by $S$, and an additive
functor $T_\thick$ induced by $T$,
making  the diagram
$$\xymatrix{(S'{\mathcal Z})/{ (S'T'{\mathcal Y})}\ar[r]^-{\pi_1}\ar[d]^S&
({\thick (S'{\mathcal Z})})/{\thick(S'T'{\mathcal Y})}\ar[d]^{S_\thick}\\
(SS'{\mathcal Z})/{( SS'T'{\mathcal Y})}\ar[r]^-{\pi_2}\ar[d]^T&
(\thick (SS'{\mathcal Z}))/{\thick( SS'T'{\mathcal Y})}\ar[d]^{T_\thick}\\
S'{\mathcal Z}/{\mathcal Y}\ar[r]^-{\pi_3}&\thick(S'{\mathcal Z})/\thick({\mathcal Y})
}$$
commutative. Moreover,
the restriction to the respective images of $\pi_1$, respectively $\pi_2$, respectively $\pi_3$
of functors $S_\thick$ and $T_\thick$ on the right
is an equivalence.
\item\label{(3)}
$S$ and $T$ induce equivalences $S_L$ and $T_L$ of additive categories fitting into the commutative diagram
$$
\xymatrix{
({\thick (S'{\mathcal Z})})_{(\thick (S'T'{\mathcal Y}))}\ar[rr]|-{S_Z}\ar@/_7pc/[ddd]_{\text{can}}&&
({\thick (SS'{\mathcal Z})})_{(\thick (SS'T'{\mathcal Y}))}\ar@/^2pc/[dddll]|-{T_Z}\\
L_{S'T'{\mathcal Y}}((S'{\mathcal Z})/{(S'T'{\mathcal Y})})\ar[r]|-{S_L}\ar[d]_{\text{can}}\ar@^{(->}[u]&
L_{SS'T'{\mathcal Y}} ((SS'{\mathcal Z})/{(SS'T'{\mathcal Y})})\ar[dl]|-{T_L}\ar@^{(->}[ur]\\
L_{{\mathcal Y}}((S'{\mathcal Z})/{\mathcal Y})\ar@^{(->}[d]\\
(\thick (S'{\mathcal Z}))_{{\mathcal Y}}&&&
}
$$
where the outer triangle consists of triangulated categories and triangle functors,
and the inner triangle are full
additive subcategories.
\item\label{(4)}
If $S$ and $T$ induce equivalences of additive categories
$$
\xymatrix{
({\thick (S'{\mathcal Z})})/{ \thick(S'T'{\mathcal Y})}\ar[rr]|-{S_\thick}\ar[d]_{\text{can}}&&
(\thick (SS'{\mathcal Z}))/{\thick( SS'T'{\mathcal Y})}\ar[dll]|-{T_\thick}\\
(\thick (S'{\mathcal Z}))/{\thick{\mathcal Y}}},
$$
then the restriction $S_Z$ of $S$ to the triangulated
category ${(\thick (S'{\mathcal Z}))}_{\thick(S'T'{\mathcal Y})}$  and
the restriction $T_Z$ of $T$ to
the triangulated category ${(\thick (SS'{\mathcal Z}))}_{\thick (SS'T'{\mathcal Y})}$
are equivalences of triangulated categories, making the diagram
$$
\xymatrix{
({\thick (S'{\mathcal Z})})_{(\thick (S'T'{\mathcal Y}))}\ar[r]|-{S_Z}\ar[d]_{\text{can}}&
({\thick (SS'{\mathcal Z})})_{(\thick (SS'T'{\mathcal Y}))}\ar[dl]|-{T_Z}\\
(\thick (S'{\mathcal Z}))_{{\mathcal Y}}}.
$$
commutative.
\end{enumerate}
\end{Theorem}

\begin{proof}
We first recall from Lemma~\ref{Zistriangulated} that $\mathcal Z$ is triangulated. By
Corollary~\ref{functorsextendtoVerdierlocalisations}
the functors coming from Theorem~\ref{Greenforadjoints} extend to functors
on the localisations. Now $S$ and $T$ are equivalences on the additive quotient
constructions. Using Proposition~\ref{functorsfrommodtoVerdier} and
Proposition~\ref{thickeningnumeratorinordinaryquotient},
the functors extend to triangle functors
$$
\xymatrix{
({\thick (S'{\mathcal Z})})_{(\thick (S'T'{\mathcal Y}))}\ar[r]|-{S_Z}\ar[d]_{\text{can}}&
({\thick (SS'{\mathcal Z})})_{(\thick (SS'T'{\mathcal Y}))}\ar[dl]|-{T_Z}\\
(\thick (S'{\mathcal Z}))_{({\mathcal Y})}}.
$$
The functors $S$ and $T$ are triangle functors on the ambient categories, and hence they induce
functors $S_\thick$ and $T_\thick$ as required.

Since Theorem~\ref{Greenforadjoints} shows that $S$ is an equivalence with
quasi-inverse $T$ on the above subcategories, the functor $T_\thick$ is also a quasi-inverse to $S_\thick$
on the images under $\pi_1$, $\pi_2$ and $\pi_3$.
Corollary~\ref{functorsextendtoVerdierlocalisations} shows item (\ref{(3)}).

Suppose now that $S$ and $T$ induce equivalences
$$
\xymatrix{
({\thick (S'{\mathcal Z})})/{( \thick(S'T'{\mathcal Y}))}\ar[rr]|-{S_\thick}\ar[d]_{\text{can}}&&
({\thick (SS'{\mathcal Z})})/{( \thick(SS'T'{\mathcal Y}))}\ar[dll]|-{T_\thick}\\
(\thick (S'{\mathcal Z}))/{\thick{\mathcal Y}}},
$$
Since by hypothesis
each object of $TSS'T'{\mathcal Y}$ is a direct factor of an object of
${\mathcal Y}$, the right vertical functor $\text{can}$ is the identity.
The restriction of the functors $S$ and $T$ in the statement of item (\ref{(4)})
are full and dense by Proposition~\ref{ordinaryquotientisogivesVerdierquotientiso}.
Since their composition $TS$ is the identity, the functors are also faithful.
The statement follows.
\end{proof}

\begin{Rem}
Note that in Theorem~\ref{Greencorrespondencefortriangulated}.(\ref{(1)}) and
\ref{Greencorrespondencefortriangulated}.(\ref{(3)}) the functor $T$ maps
from the localisation at the thick subcategory of images under $S$ to
the localisation at the thick subcategory of images under $T$. Note that by
Proposition~\ref{HigmansLemmaForAdjointFunctors} we get a functor from the localisation at the thick
subcategory generated by $T$-relative injective objects to the localisation at the thick
subcategory generated by $S$-relative projective objects.
\end{Rem}

\begin{Rem}
Consider the special situation when $G$ is a finite group and $k$ is a field of
characteristic $p>0$.
Then, following Carlson, Peng and Wheeler~\cite{Carlson-Peng-Wheeler} the classical
Green correspondence is an equivalence of full additive subcategories of triangulated categories.

More precisely, let $D$ be a $p$-subgroup of $G$ and let $H$ be a subgroup of $G$ containing
$N_G(D)$, the normalizer of $D$ in $G$.

Consider ${\mathcal G}=kG-\ul{mod}$, ${\mathcal H}=kN_G(D)-\ul{mod}$, and
${\mathcal D}=kD-\ul{mod}$, the stable categories of $kG$-modules, $kH$-modules,
and $kD$-modules.
Here the stable categories are taken modulo morphisms factoring through projective modules.
Let $$S=kG\otimes_{kH}-=\ind_H^G:kH-\ul{mod}\lra kG-\ul{mod}$$ and
$$S'=kH\otimes_{kD}-=\ind_D^H:kD-\ul{mod}\lra kH-\ul{mod}$$
be the induction functors.
These have left and right adjoints, namely the restriction
$$T:=Hom_{kH}(kG,-)=\res^G_H:kG-\ul{mod}\lra kH-\ul{mod}$$ is left and right adjoint to $S$.
Similarly,
$$T':=Hom_{kH}(kG,-)=\res^G_H:kH-\ul{mod}\lra kD-\ul{mod}$$ is left and right adjoint to $S'$.

Since group algebras are symmetric, following Remark~\ref{Frobeniusforleftequalsrightadjoint}
the stable categories
${\mathcal G}=kG-\ul{mod}$, ${\mathcal H}=kN_G(D)-\ul{mod}$, and
${\mathcal D}=kD-\ul{mod}$ are triangulated and moreover, the functors
$\ind_H^G, \ind_D^H,\res^G_H,\res^H_D$ come from exact functors of the corresponding
module categories, and hence are triangle functors. Further,
$$\res^G_H\ind^G_H=\id_{kH-mod}\oplus U$$
for $U=\bigoplus_{HgH\in H\backslash G/H\setminus\{H\}}{} kHgH\otimes_{kH}-$.

Therefore Theorem~\ref{Greencorrespondencefortriangulated} applies for
appropriate choices of $\mathcal Y$. Following \cite[page 311, Section 3]{AuslanderKleiner}
we fix a collection $\mathfrak Y$ of subgroups of $H$, closed under
$H$-conjugation and under taking subgroups, we consider
${\mathcal Y}$ the full subcategory of ${\mathcal H}$ given by
$\textup{ind}_{\mathfrak Y}^H$, i.e. those $kH$-modules induced from $kY$-modules for
some $Y\in{\mathfrak Y}$. By \cite[Corollary 3.4 (a) and (b)]{AuslanderKleiner}
we may put ${\mathfrak Y}:=\{Y\;|\;Y\leq H\cap gDg^{-1}\;;\;g\in G\setminus H\}$
which satisfies the hypotheses of Theorem~\ref{Greencorrespondencefortriangulated}.
Moreover, the functors $S_L$ and $T_L$ in item $(3)$ of
Theorem~\ref{Greencorrespondencefortriangulated} implies the classical Green correspondence.
Furthermore, the bijection of indecomposable $kG$-modules and $kH$-modules
with vertex $D$ is the restriction of a triangle functor between triangulated categories,
namely the Verdier localisation of triangulated subcategories.

However, if $D$ is TI, i.e. $D\cap D^g\in\{1,D\}$ for all $g\in G$,
the stable categories involved in the theorem are the usual stable categories modulo
projectives, which are already triangulated, and by the universal
property of the Verdier localisation (\cite[$\S 2$, no 3]{SGA412} or
\cite[Chapitre II, Corollaire 2.2.11.c]{Verdier})
there is an inverse functor to $L$ (which was introduced in
Lemma~\ref{naturalquotientinduces}).

By the same argument, for general $D$, the Verdier localisation in item $(3)$ of
Theorem~\ref{Greencorrespondencefortriangulated} is the $W$-stable category from
Carlson-Peng-Wheeler~\cite{Carlson-Peng-Wheeler} (cf also Grime \cite[Example 3.6]{Grime}).

\end{Rem}

\section{Tensor triangulated categories---Green correspondence abstractly and for group rings}

\label{tensortriangulatedsection}

We had to deal with thick subcategories of triangulated categories. Our main model was the case of
versions of derived or stable categories of group rings. Classification results are known in
this case, but mainly in presence of an additional monoidal structure.

\subsection{Recall Balmer's results}

We first recall some results from Balmer \cite{Balmer}.
\begin{itemize}
\item A tensor triangulated category $\mathcal K$ is an
\begin{itemize}
\item essentially small
\item
triangulated category $\mathcal K$
together with a
\item symmetric monoidal structure $({\mathcal K},\otimes,1)$,
\item such that the functor $\otimes :{\mathcal K}\times {\mathcal K}\ra {\mathcal K}$
is assumed to be exact in each variable.
\end{itemize}
\item A tensor triangulated functor is an exact functor between tensor triangulated categories
sending the identity object to the identity object and respecting the monoidal structures.
\item
A $\otimes$-ideal $\mathcal P$ of ${\mathcal K}$ is a
\begin{itemize}
\item thick triangulated
subcategory
\item such that if an object $M$ is in $\mathcal P$ and $X$ is an object in
$\mathcal K$, then $M\otimes X$ is in $\mathcal P$.
\end{itemize}
\item An ideal $\mathcal P$ is prime if
$A\otimes B$ being an object in $\mathcal P$ if and only if $A$ is an object in $\mathcal P$
or $B$ is an object in $\mathcal P$.
\item The spectrum $\Spec({\mathcal K})$ is defined to be the set (!) of prime ideals
of $\mathcal K$.
\item The support of an object $M$ of $\mathcal K$ is
$$\supp_{\mathcal K}(M):=\{{\mathcal P}\in \Spec({\mathcal K})\;|\;
M\textup{ is not an object of }{\mathcal P}\}.$$
\item For any family of objects ${\mathcal S}$ of $\mathcal K$ let
$Z({\mathcal S}):=\{{\mathcal P}\in \Spec({\mathcal K})\;|\;{\mathcal S}\cap{\mathcal P}=\emptyset\}$.
The sets $Z({\mathcal S})$ form the closed sets of a topology, the Zariski topology on $\Spec({\mathcal K})$.
\item The radical $\sqrt{\mathcal P}$ of an ideal $\mathcal P$ is the class
of objects $M$ in $\mathcal K$
such that there is $n\in\N$ so that $M^{\otimes n}$ is an object of $\mathcal P$. An ideal $\mathcal P$ is called radical if $\sqrt{\mathcal P}=\mathcal P$.
\end{itemize}

One of the main results of \cite{Balmer} is

\begin{Theorem}\cite[Theorem 4.10]{Balmer} Let ${\mathfrak S}({\mathcal K})$ denote the
subsets $Y\subseteq \Spec({\mathcal K})$ such that $Y=\bigcup_{i\in I}Y_i$ with all $Y_i$
closed and $\Spec({\mathcal K})\setminus Y_i$ quasicompact. Let ${\mathfrak R}({\mathcal K})$
be the set of radical thick $\otimes$-ideals of $\mathcal K$.
Then the following maps are mutually inverse bijections
\begin{eqnarray*}
{\mathfrak S}({\mathcal K})&\longleftrightarrow& {\mathfrak R}({\mathcal K})\\
Y&\mapsto&{\mathcal K}_Y:=\{M\in{\mathcal K}\;|\;\supp(M)\subseteq Y\}\\
\bigcup_{M\in{\mathcal J}}\supp(M)=:\supp({\mathcal J})&\leftmapsto&{\mathcal J}
\end{eqnarray*}
\end{Theorem}

\subsection{Green correspondence of the spectrum in a tensor triangulated category}

Recall that, following \cite{Etingofetal}
a tensor subcategory of a tensor category still has a unit element.
We shall need to define a concept without this restriction since for our natural examples
we do not necessarily have a unit element.
Note that a semigroup is a set with a binary associative structure, and a monoid is a
semigroup with a unit. We transport this vocabulary to the world of tensor categories
under the name of semigroup category (cf \cite{Borychenko}).

\begin{Def}\label{semigroupcategorydef}
\begin{itemize}
\item A semigroup category is a category $\mathcal C$ with a symmetric binary operation
$\otimes:{\mathcal C}\times{\mathcal C}\ra{\mathcal C}$ satisfying the associative pentagon axiom.
\item
A triangulated semigroup category is
\begin{itemize}
\item an essentially small triangulated category ${\mathcal C}$,
\item which is in addition a semigroup category $({\mathcal C},\otimes)$
\item such that $\otimes:{\mathcal C}\times{\mathcal C}\ra{\mathcal C}$ is exact in each variable.
    \end{itemize}
\item
A $\otimes$-ideal $\mathcal P$ of a triangulated semigroup category ${\mathcal C}$ is a
\begin{itemize}
\item thick triangulated subcategory
\item such that if an object $M$ is in $\mathcal P$ and $X$ is an object in
$\mathcal C$, then $M\otimes X$ is in $\mathcal P$.
\end{itemize}
\item
Let $({\mathcal C},\otimes)$ and $({\mathcal D},\otimes)$ be semigroup categories.
\begin{itemize}
\item A functor
$F:{\mathcal C}\ra{\mathcal D}$ is called semi-tensor functor if $F$ allows a natural equivalence
$J:F(V\otimes W)\ra F(V)\otimes F(W)$ which satisfies
the associahedron diagram \cite[Diagram 2.23]{Etingofetal}.
\item If $({\mathcal C},\otimes)$ and $({\mathcal D},\otimes)$ are triangulated
semigroup categories, then a semi-tensor functor $F:{\mathcal C}\ra{\mathcal D}$ is called
    triangle semi-tensor functor if $F$ is in addition a triangle functor.
\end{itemize}
\end{itemize}
\end{Def}

\begin{Lemma}\label{tensoridealinverseimage}
If ${\mathcal Y}$ is a thick semigroup triangulated subcategory (respectively
$\otimes$-ideal) of a triangulated semigroup category $\mathcal H$,
and if $F:{\mathcal D}\lra{\mathcal H}$ is a triangle semi-tensor functor, then
${\mathcal Z}:=F^{-1}({\mathcal Y})$ is again a thick semigroup triangulated subcategory
(resp. $\otimes$-ideal) of $\mathcal D$.
\end{Lemma}

\begin{proof} By Proposition~\ref{Zistriangulated} $\mathcal Z$ is triangulated.
By definition, $\mathcal Z$ is thick. We need to show that $\mathcal Z$ is a semigroup tensor category.
Let $X$ be an object of ${\mathcal Z}$ and $Y$ be an object of ${\mathcal Z}$
(resp. $\mathcal H$). Then there are objects $X'$ and $Y'$ such that
$F(X)\oplus X'$ and $F(Y)\oplus Y'$ are objects of ${\mathcal Y}$.
But then
\begin{eqnarray*}
(F(X)\oplus X')\otimes (F(Y)\oplus Y')&\simeq&(F(X)\otimes F(Y))\\
&&\oplus (X'\otimes F(Y))\oplus (F(X)\otimes Y')\oplus (X'\otimes Y')\\
&\simeq&F(X\otimes Y)\\
&&\oplus (X'\otimes F(Y))
\oplus (F(X)\otimes Y')\oplus (X'\otimes Y')
\end{eqnarray*}
Since $\mathcal Y$ is tensor triangulated, $F(X\otimes Y)$ is a direct factor of
the object $(F(X)\oplus X')\otimes (F(Y)\oplus Y')$ of $\mathcal Y$. Therefore
$X\otimes Y$ is an object of $\mathcal Z$.
\end{proof}

\begin{Lemma}\label{localisationonidealistensortriangulated}
If ${\mathcal H}$ is a tensor triangulated category, if ${\mathcal Y}$ is a $\otimes$-ideal
in ${\mathcal H}$, then the tensor triangulated structure on $\mathcal H$
induces a  tensor triangulated structure on ${\mathcal H}_{\mathcal Y}$.
Moreover, the natural functor $\nu:{\mathcal H}\ra {\mathcal H}_{\mathcal Y}$
is a functor of tensor triangulated categories.
\end{Lemma}

\begin{proof}
The objects of ${\mathcal H}$ coincides with the objects of ${\mathcal H}_{\mathcal Y}$.
We need to define a tensor product $\ol\otimes$ on ${\mathcal H}_{\mathcal Y}$.
Denote by $\nu:{\mathcal H}\ra {\mathcal H}_{\mathcal Y}$ the natural functor.
We define for any two objects $M,N$  in ${\mathcal H}_{\mathcal Y}$ the object
$M\ol{\otimes} N:=\nu(M\otimes N)$ in
${\mathcal H}_{\mathcal Y}$.

Since $\mathcal Y$ is an ideal, this construction is also well-defined on morphisms.
Since $\id_{\mathcal H}$ is the neutral element of $\otimes$, we
get $\nu(\id_{\mathcal H})$ is the
neutral element of $\ol\otimes$.
Since $\otimes$ is monoidal symmetric, also $\ol\otimes$ is monoidal symmetric.
The functor is tensor triangulated by construction.
\end{proof}

Recall that every thick subcategory is a full
triangulated subcategory, but
a full triangulated subcategory is thick only if it is in addition
closed under taking direct summands.
A full triangulated subcategory $\mathcal A$ of $\mathcal D$ is
strict if any object of $\mathcal D$
which is isomorphic in $\mathcal D$ to an object in $\mathcal A$
is also an object of $\mathcal A$.

\begin{Prop}
Let $({\mathcal T},\otimes,1)$ be a tensor triangulated category and let
$\mathcal P$ and $\mathcal Q$ be $\otimes$-ideals of $\mathcal T$. Suppose moreover
that ${\mathcal Q}$ is a full triangulated subcategory of $\mathcal P$. Suppose that $\mathcal P$
is strictly full in $\mathcal T$.

Then the following hold.
\begin{itemize}
\item The tensor triangulated structure on $\mathcal T$ induces a
tensor triangulated structure
$\ol\otimes$ on the Verdier localisation ${\mathcal T}_{\mathcal Q}$.
\item
Furthermore, consider the natural functor $\nu:{\mathcal T}\lra {\mathcal T}_{\mathcal Q}$.
Let $\nu'$ be the restriction of $\nu$ to ${\mathcal P}$, as indicated in the commutative diagram
$$\xymatrix{{\mathcal T}\ar[r]^\nu &{\mathcal T}_{\mathcal Q}\\
{\mathcal P}\ar@{^{(}->}[u]\ar[r]^{\nu'}&\nu({\mathcal P})\ar@{^{(}->}[u]}$$
Denote by ${\mathcal P}_{({\mathcal Q})}$ the image $\nu({\mathcal P})$ of $\mathcal P$
in ${\mathcal T}_{\mathcal Q}$
under $\nu$, and denote by ${\mathcal P}_{{\mathcal Q}}$ the
Verdier localisation of $\mathcal P$ at $\mathcal Q$. Then
$${\mathcal P}_{({\mathcal Q})}={\mathcal P}_{\mathcal Q}.$$
\item
${\mathcal P}_{\mathcal Q}$ is a $\ol\otimes$-ideal of ${\mathcal T}_{\mathcal Q}$.
\end{itemize}
\end{Prop}

\begin{proof}
Lemma~\ref{localisationonidealistensortriangulated} is precisely the first statement.

\medskip

Denote by $\iota:{\mathcal P}\ra{\mathcal T}$ the inclusion functor.
As for the second statement we have the Verdier localisation  ${\mathcal P}_{{\mathcal Q}}$
of $\mathcal P$ at $\mathcal Q$. Denote by $\mu:{\mathcal P}\ra{\mathcal P}_{{\mathcal Q}}$
the natural functor. Then, the universal property of Verdier localisations (\cite[$\S 2$, no 3]{SGA412} or
\cite[Chapitre II, Corollaire 2.2.11.c]{Verdier}) induces a unique
functor $\sigma:{\mathcal P}_{{\mathcal Q}}\ra{\mathcal P}_{({\mathcal Q})}$
such that $\sigma\circ\mu=\nu\circ\iota$. This shows that the functor $\sigma$ is dense
since $\mu$, $\nu$, $\iota$ are the identity on objects.

We need to show that $\sigma$ is fully faithful.
Let $Z$ be an object of $\mathcal T$, let $P_1$ and $P_2$ be objects of $\mathcal P$,
and a diagram of morphisms of $\mathcal T$
$$
\xymatrix{&Z\ar[dr]^{\alpha}\ar[dl]_{\gamma}\\P_1&&P_2}
$$
representing a morphism $\omega$ in ${\mathcal T}_{\mathcal Q}(P_1,P_2)$.
If $\gamma$ has cone in $\mathcal Q$, since $\mathcal P$ is triangulated, and
since $\mathcal Q$ is a triangulated subcategory of $\mathcal P$,
also $Z$ is isomorphic to an object of $\mathcal P$, and since $\mathcal P$ is
strictly full in $\mathcal T$,
the object $Z$ is actually an object of $\mathcal P$. Hence
$\sigma$ is full.

If $\lambda$ is represented by
$$
\xymatrix{&Z\ar[dr]^{\alpha}\ar[dl]_{\gamma}\\P_1&&P_2}
$$
for some object $Z$ of $\mathcal P$,
and if $\sigma(\lambda)=0$ in ${\mathcal P}_{\mathcal Q}$, then there is an object $Z'$ of $\mathcal T$
and a morphism $\delta:Z'\ra Z$ with $\cone(\delta)$ in $\mathcal Q$, and with $\alpha\circ\delta=0$.
But, again $\mathcal Q$ is a triangulated subcategory of $\mathcal P$, and $\mathcal P$
being strictly full triangulated subcategory of $\mathcal T$
implies $Z'$ is an object of $\mathcal P$. Since $\mathcal P$ is a full
subcategory of $\mathcal T$, the morphism $\delta$ is actually already in $\mathcal P$. Hence $\lambda=0$.
This shows that $\sigma$ is faithful. Altogether we get the second statement.

\medskip

Since $\mathcal P$ is a $\otimes$-ideal, for any $P$ in $\mathcal P$, and any $X$ in $\mathcal T$ we get
$P\otimes X$ is in $\mathcal P$. Hence $$\nu(P)\ol\otimes\nu(X)=\nu(P\otimes X)$$ is an object of
${\mathcal P}_{({\mathcal Q})}={\mathcal P}_{\mathcal Q}$. This proves the third statement.
\end{proof}

\subsection{Thick tensor triangulated categories and tensor ideals in the special case of group rings }

\label{specialcaseofgrouprings}

Various results are known for classification of thick subcategories of various
triangulated categories (cf e.g.
\cite{Thomason, Benson-Carlson-Rickard, Benson-Iyengar-Krause, Carlson-Iyengar,
FriedlanderPevtsova}),
giving mostly a parametrisation with certain subsets of support varieties.
For a fixed, essentially small triangulated category $\mathcal D$
a general result describing the relation between full triangulated essentially small
subcategories $\mathcal A$ and $\thick({\mathcal A})={\mathcal D}$
is given by Thomason.

\begin{Theorem} (Thomason \cite[Theorem 2.1]{Thomason})
Let $\mathcal D$ be an essentially small triangulated category.
Consider the set $\mathfrak U$ of
strictly full triangulated subcategories
$\mathcal A$ in $\mathcal D$, having the property that each
object of $\mathcal D$ is isomorphic to a
direct summand of an object in $\mathcal A$. Then $\mathfrak U$
is in bijection with
the set of subgroups of the  Grothendieck group $K_0({\mathcal D})$.
The isomorphism is given by mapping $\mathcal A$ to the
subgroup $K_0({\mathcal A})$
of $K_0({\mathcal D})$.
\end{Theorem}

Thomason also gives \cite[Theorem 3.15]{Thomason} a classification of
tensor triangulated thick subcategories of the derived category of perfect
complexes over a quasi-compact quasi-separated scheme.

We focus on those dealing with group rings.
Let $k$ be an algebraically closed field of characteristic $p>0$ and let $G$ be a finite
group with order divisible by $p$.
Let $H^\bullet(G)$ be $\bigoplus_{i\geq 0}H^{2i}(G,k)$ if $p$ is odd, and
$H^\bullet(G)=H^*(G,k)$ if $p=2$. Then $H^\bullet(G)$ is a graded commutative algebra,
and $Ext_{kG}^*(M,M)$ is a finitely generated $H^\bullet(G)$-module.
Let $V_G(k)$ be the maximal ideal spectrum of
$H^\bullet(G)$.
A set $\mathcal X$ of closed subvarieties of $V_G(k)$ is said to be closed under specialisation if
whenever
$W\in{\mathcal X}$ and $W'\subseteq W$, then we also get $W'\in{\mathcal X}$.
For a set $\mathcal X$ of closed subvarieties of $V_G(k)$ which is closed under specialisation
we let ${\mathcal C}({\mathcal X})$ be the thick subcategory of $kG-\ul{mod}$
consisting of modules $M$ with
$$V_G(M):=\{\mathfrak m\in V_G(k)\;|\;\textup{Ann}_{H^\bullet(G)}(Ext^*_{kG}(M,M))
\subseteq{\mathfrak m}\}\in{\mathcal X}.$$

Benson, Carlson and Rickard showed in \cite{Benson-Carlson-Rickard} the following.

\begin{Theorem} \cite[Theorem 3.4]{Benson-Carlson-Rickard}
Let $k$ be an algebraically closed field, and let $G$ be a finite group.
Let $V_G(k)$ be the maximal ideal spectrum of $H^\bullet(G)$.
Then the thick tensor ideals $I$ in $kG-\ul{mod}$ are of the form ${\mathcal C}({\mathcal X})$
for some non empty set $\mathcal X$ of homogeneous subvarieties of $V_G(k)$, closed under
specialisation and finite unions.
\end{Theorem}

Carlson and Iyengar \cite{Carlson-Iyengar} determined the thick subcategories of
the derived category of the group algebra of a finite group.
For each object $M$ of
$D^b(kG)$ there is a morphism of $k$-algebras $H^*(G,k)\ra Ext_{kG}^*(M,M)$.
Again $Ext_{kG}^*(M,M)$ is a finitely generated $H^\bullet(G)$-module.
Then
$$V_{D^b(kG)}(M):=\textup{Supp}_{H^\bullet(G)}(M):=\{\wp\in
\textup{Spec}(H^\bullet(G))\;|\;H(M_\wp)\neq 0\}
\subseteq \textup{Spec}(H^\bullet(G)).$$

\begin{Theorem}\cite[Theorem 6.6 and Corollary 6.7]{Carlson-Iyengar}
For an algebraically closed field $k$ of characteristic $p>0$ and a finite
group $G$ with order divisible by $p$, and two objects $M$ and $N$ in $D^b(kG)$ with
$V_{D^b(kG)}(M)\subseteq V_{D^b(kG)}(N)$ we have that $M$ is in the thick tensor ideal
generated by $N$.
In particular, if $\mathbf C$ is a thick tensor ideal of $D^b(kG)$, then there is a
specialisation closed subset $V$ of $V_{D^b(kG)}(k)$ such that $\mathbf C$
equals the subcategory obtained by all those $M$ in $D^b(kG)$ with $V_{D^b(kG)}(M)\subseteq V$.
\end{Theorem}

Carlson \cite{Carlson} studied thick subcategories of what he calls
relatively stable categories of group rings.
Let $\mathfrak H$ be a set of subgroups of $G$. A $kG$-module $M$ is called $\mathfrak H$-projective
if $M$ is $\dar^G_H$-relative projective for all $H\in{\mathfrak H}$.
It is classical that a module $M$ is $\mathfrak H$-projective if and only if
$M$ is a direct summand of modules which are induced from modules over elements of $\mathfrak H$.
The category $kG-\ul{mod}_{\mathfrak H}$ has the same objects as $kG-mod$.
However, the set of morphisms from $M$ to $N$ is the set of equivalence classes of
$kG$-module morphisms modulo those factoring through $\mathfrak H$-projective modules.

Carlson, Peng and Wheeler \cite[Theorem 6.2]{Carlson-Peng-Wheeler} show that
$kG-\ul{mod}_{\mathfrak H}$ is actually a triangulated category. Moreover,
an immediate consequence is that Green correspondence is the restriction of a functor
between triangulated categories to certain subcategories, namely those full subcategories
generated by indecomposable modules with a specific vertex. This restriction is then
an equivalence of categories.

Benson and Wheeler extend the concept to infinitely generated modules, and
show in \cite[Proposition 2.3]{Benson-Wheeler} that we get again triangulated categories and
a Green correspondence, which is an equivalence between triangulated categories.

These results are special cases of our more general approach, when applied to
bounded derived categories of modules over the respective group rings and
appropriate choices for $\mathcal Y$, as is shown by the following.

\begin{Prop}\label{relativelystableislocalisation}
Let $D$ be a $p$-subgroup of $G$, let $H$ be a subgroup of $G$ containing $N_G(D)$,
and let ${\mathfrak Y}:=\{S\leq H\cap {}^gD\;|\;g\in G\setminus H\}$ as well as
${\mathfrak X}:=\{S\leq D\cap {}^gD\;|\;g\in G\setminus H\}$.

Then for $\mathcal Y$ being the class of complexes having indecomposable factors with
vertex in $\mathfrak Y$, the natural functors
$$L_{\mathcal SS'T'Y}:kG-\ul{mod}_{\mathfrak X}\ra D^b(kG)_{\thick(SS'T'{\mathcal Y})}$$
and
$$L_{\mathcal Y}:kH-\ul{mod}_{\mathfrak Y}\ra D^b(kH)_{\thick({\mathcal Y})}$$
are equivalences of triangulated categories.
\end{Prop}

\begin{proof} Using \cite[Lemma 4.1]{Wang-Zhang},
Lemma~\ref{naturalquotientinduces} defines $L_{SS'T' \mathcal Y}$ and
$L_{\mathcal Y}$.
Since the subcategory of bounded complexes of finitely generated projectives is a
subcategory of $\mathcal Y$, the category $D^b(kG)_{\thick(SS'T'{\mathcal Y})}$
is a localisation of the singularity category $D_{sg}(kG)$ of $kG$. The singularity category
of a self-injective algebra is just the stable category of the algebra modulo
projective-injectives (cf
Keller-Vossieck~\cite{KellerVossieck}, Rickard~\cite{Rickardstable},
Buchweitz~\cite{Buchweitz}). Likewise $D^b(kH)_{\thick({\mathcal Y})}$ is a
localisation of the stable category of $kH$.
Since the categories $kG-\ul{mod}_{\mathfrak X}$ and
$kH-\ul{mod}_{\mathfrak Y}$ are triangulated, the universal property of
the Verdier localisation (\cite[$\S 2$, no 3]{SGA412} or
\cite[Chapitre II, Corollaire 2.2.11.c]{Verdier}) gives the quasi-inverse functors
to $L_{SS'T' \mathcal Y}$ and $L_{\mathcal Y}$ respectively.
\end{proof}

Observe that now Theorem~\ref{Greencorrespondencefortriangulated} item (3)
gives Carlson, Peng and Wheeler's theorem.

Moreover, Harris in \cite{Harris} and independently Wang and Zhang  in \cite{Wang-Zhang}
give a blockwise version of the Green correspondence.

Localising subcategories are a vast variety in this setting.
Carlson \cite{Carlson}
shows for example that for $p=2$ and a collection $\mathcal C$ of subgroups $H$ of $G$
all of which containing an elementary abelian subgroup of rank at least $2$, then the
spectrum of the relatively
$\mathcal C$-stable category is not Noetherian.

Note that for a non principal block we do not get a monoidal category but only
a semigroup category in the sense of Definition~\ref{semigroupcategorydef}.
Indeed, the unit element is the trivial module, which belongs to the principal block.


\begin{thebibliography}{88}

\bibitem{AuslanderKleiner}
Maurice Auslander and Mark Kleiner, {\em Adjoint functors and an Extension of the Green
Correspondence for Group Representations}, Advances in Mathematics {\bf 104} (1994) 297-314.

\bibitem{AR}
Maurice Auslander and Idun Reiten, {\em Homologically finite subcategories}, In H. Tachikawa and S. Brenner (Eds.), Representations of Algebras and Related Topics (London Mathematical Society Lecture Note Series {\bf 168}, 1-42)

\bibitem{Balmer}
Paul Balmer, {\em  The spectrum of prime ideals in tensor triangulated categories,}
Journal f\"ur die reine und angewandte Mathematik {\bf 588} (2005), 149-168.

\bibitem{BBD} Alexander A. Beilinson, Joseph Bernstein, Pierre Deligne, {\em Faisceaux pervers}.
Ast\'erisque 100, vol. 1. Soci\'et\'e Math\'ematique de France (1982).

\bibitem{BeligiannisMarmaridis}
Apostolidis Beligiannis and Nikolaos Marmaridis, {\em Left triangulated categories
arising from contravariantly finite subcategories},
Communications in Algebra {\bf 22} (12), (1994) 5021-5036.

\bibitem{Benson-Wheeler}
David J. Benson and Wayne Wheeler, {\em The Green correspondence for infinitely generated modules},
Journal of the London Mathematical Society {\bf 63} (2001) 69-82.

\bibitem{Benson-Carlson-Rickard}
David J. Benson, Jon F. Carlson and Jeremy Rickard, {\em Thick subcategories of the stable
module category}, Fundamenta Mathematicae {\bf 153} (1997) 59-80.

\bibitem{Benson-Iyengar-Krause}
David J. Benson, Srikanth Iyengar and Henning Krause, {\em Stratifying modular
representations of finite groups}, Annals of Mathematics {\bf 174} (2011) 1643-1684.

\bibitem{Borychenko}
Mitya Boyarchenko, {\em Associative constraints in monoidal categories},
unpublished manuscript, University of Chicago 2011.

\bibitem{BroueHigman}
Michel Brou\'e, {\em Higman's criterion revisited}, Michigan Mathematical Journal {\bf 58} (2009), no.1, 125-179.

\bibitem{Buchweitz}
Ragnar Buchweitz, {\em Maximal Cohen-Macaulay modules and Tate-cohomology over Gorenstein rings},
unpublished manuscript, Hannover 1987.

\bibitem{Buehler}
Theo B\"uhler, {\em Exact categories}, Expositiones Math. {\bf 28} (2010) 1-69.

\bibitem{Carlson-Peng-Wheeler}
Jon F. Carlson, Chuang Peng, and Wayne Wheeler, {\em Transfer maps and Virtual Projectivity},
Journal of Algebra {\bf 204} (1998) 286-311.

\bibitem{Carlson}
Jon F. Carlson, {\em Thick subcategories of the relative stable category}, in
'Geometric and Topological Aspects of the Representation Theory of Finite Groups,
PIMS Summer School and Workshop, July 27-August 5, 2016; Springer Proceedings in
Mathematics and Statistics 2018, pages 25-49.

\bibitem{Carlson-Iyengar}
Jon F. Carlson and Srikanth Iyengar, {\em Thick subcategories of the bounded derived
category of a finite group}, Transactions of the
American Mathematical Society {\bf 367} (2015) 2703-2717.

\bibitem{Eilenberg-Moore}
Samuel Eilenberg and J. C. Moore, {\em Foundations of relative homological algebra}
Memoirs of the American Mathematical Society {\bf 55} (1966).

\bibitem{Etingofetal}
Pavel Etingof, Shlomo Gelaki, Dmitri Nikshych, and Victor Ostrik, {\sc Tensor Categories},
American Mathematical Society Surveys and Monographs 205, Providence Rhode Island 2015.

\bibitem{FriedlanderPevtsova}
Eric Friedlander and Julia Pevtsova, {\em $\Pi$-supports for finite group schemes over a field},
Duke Mathematical Journal {\bf 139} (2007) 317-368.

\bibitem{Green}
James Alexander Green, {\em On the indecomposable representations of a finite group},
Mathematische Zeitschrift {\bf 70} (1959) 430-445.

\bibitem{Grime}
Matthew Grime, {\em Adjoint functors and triangulated categories},
Communications in Algebra {\bf 36} (2008) 3589-3607.

\bibitem{Happel}
Dieter Happel, {\em Triangulated categories in the representation theory of finite
dimensional algebras}, London Mathematical Society Lecture Note Series {\bf 119} (1988).


\bibitem{Harris}
Morton E. Harris, {\em A block extension of categorical results in the Green correspondence},
Journal of group theory {\bf 17} (2014) 1117-1131.

\bibitem{Hochschild}
Gerhard Hochschild, {\em Relative Homological Algebra}, Transactions of the American Mathematical Society
{\bf 82} (1956) 246-269.


\bibitem{KellerVossieck}
Bernhard Keller and Dieter Vossieck, {\em Sous les cat\'egories d\'eriv\'ees}, Comptes Rendus de l'Acad\'emie des Sciences Paris {\bf 305} (1987) 225-228.

\bibitem{Maclane}
Saunders Maclane, {\sc Categories for the Working Mathematician},
Second edition, Springer Verlag 1978.

\bibitem{RickardMoritatheorem}
Jeremy Rickard, {\em Morita theory for derived categories}, Journal of the 
London Mathematical Society (2) {\bf 39} (1989) 436-456.

\bibitem{Rickardstable}
Jeremy Rickard, {\em Derived categories and stable equivalences},
Journal of Pure and Applied Algebra {\bf 61} (1989) 303-317.

\bibitem{Objective}
Claus M. Ringel and Pu Zhang, {\em Objective triangle functors},
Science China Mathematics, February 2015 Vol. 58, no2; 221-232

\bibitem{Stacks}
The Stacks project authors, {\em The Stacks Project}, {\tt{https://stacks.math.columbia.edu}}, year: 2019.

\bibitem{Thomason}
Robert W. Thomason, {\em The classification of triangulated subcategories},
Compositio Mathematica {\bf 105} (1997) 1-27.

\bibitem{SGA412}
Jean-Louis Verdier, {\em Cat\'egories d\'eriv\'ees: Quelques r\'esultats (\'etat 0)},
in SGA $4\frac12$-Cohomologie \'etale (Pierre Deligne) Springer Lecture notes {\bf 569} (1977) 262-311.

\bibitem{Verdier}
Jean-Louis Verdier, {\sc Des cat\'egories d\'eriv\'ees des cat\'egories ab\'eliennes},
Ast\'erisque {\bf 239} (1996), Soci\'et\'e Math\'ematique de France.

\bibitem{Wang-Zhang}
Lizhong Wang and Jiping Zhang, {\em Relatively stable equivalences
of Morita type for blocks}, Journal of pure and applied Algebra {\bf 222} (2018) no 9, 2703-2717.

\bibitem{reptheobuch}
Alexander Zimmermann, {\sc Representation theory: a homological algebra point of view},
Springer Verlag Cham, 2014.

\bibitem{RemarksOnGreen}
Alexander Zimmermann, {\em Remarks on a triangulated version of Auslander-Kleiner's Green correspondence}, preprint 2020.

\end{thebibliography}
\end{document}